\input amstex.tex

\input amsppt.sty

\TagsAsMath

\magnification=1200

\hsize=5.0in\vsize=7.0in

\hoffset=0.2in\voffset=0cm

\nonstopmode

\document

\input amstex.tex
\input amsppt.sty
\TagsAsMath \NoRunningHeads \magnification=1200
\hsize=5.0in\vsize=7.0in \hoffset=0.2in\voffset=0cm \nonstopmode

\document

\topmatter

\title{Orbitally but not asymptotically stable  ground states for the discrete
 NLS}
\endtitle

\author
Scipio Cuccagna
\endauthor

\address
DISMI University of Modena and Reggio Emilia, via Amendola 2,
Padiglione Morselli, Reggio Emilia 42100 Italy\endaddress \email
cuccagna.scipio\@unimore.it \endemail

\abstract We  consider examples of discrete nonlinear Schr\"odinger
equations  in $\Bbb Z$ admitting ground states which are orbitally
but not asymptotically stable in $\ell ^2(\Bbb Z)$. The ground
states contain internal modes  which decouple from the continuous
modes. The absence of leaking of energy from discrete to continues
modes leads to an almost conservation and perpetual oscillation of
the discrete modes. This is quite different from what is known for
nonlinear Schr\"odinger equations in $\Bbb R^d$. We do not
investigate connections with work on quasi periodic solutions as in
\cite{JA,BV}.

\endabstract

\endtopmatter

 \head \S 1 Introduction \endhead

 We
consider the   discrete Laplacian $\Delta$ in $\Bbb Z$ defined by
$$(\Delta u)(n)= u(n+1)+u(n-1)-2u(n) .$$
 In $\ell ^2 (\Bbb
Z)$ we have for the spectrum $\sigma (-\Delta
 )=[0,4 ]$. Let for $\langle n \rangle=\sqrt{1+n^2}$
$$\aligned & \ell ^{p,\sigma}(\Bbb Z)=\{ u=\{ u_n\} :
\| u\|  _{\ell ^{p,\sigma}} ^p=\sum _{n\in \Bbb Z} \langle n \rangle
^{p\sigma }|u(n)|^p<\infty \} \text{ for $p\in [1,\infty )$}\\& \ell
^{\infty ,\sigma}(\Bbb Z)=\{ u=\{ u(n)\} : \| u\|  _{\ell
^{\infty,\sigma}}
 =\sup _{n\in \Bbb Z} \langle n \rangle ^{ \sigma }|u(n)| <\infty
\}  .\endaligned $$ We set $\ell ^{p,\sigma}=\ell ^{p,\sigma}(\Bbb
Z)$ and $\ell ^{p  } =\ell ^{p,0}.$ We consider a potential $q=\{
q(n), n\in \Bbb Z \}$ with $q(n)\in \Bbb R$ for all $n$.
 We consider the discrete Schr\"odinger operator $H$
$$(Hu)(n)=-( \Delta  u)(n)+q(n)u(n).\tag 1.1$$
We   assume:

{\item{(H1)}} $|q (n)|\le C e^{-|n|} $.

{\item{(H2)}} The points 0 and 4 are not resonances of $H$.

{\item{(H3)}} $\sigma _d(H)$ consists of exactly two eigenvalues
$-E_0<0$ and $E_1>4$.

  Here $\lambda =0$ (resp.4) is a resonance if there
is a nontrivial $u\in \ell ^\infty $ with $Hu=\lambda u$. We will
see in Appendix A that there are operators satisfying (H1)--(H3). By
Lemma 5.3 \cite{CT},   $\dim \ker (H+E_0)= 1 $ in $\ell ^2$. We
denote by $\varphi _0(n)$ a generator of $\ker (H+E_0) $ normalized
so that $\| \varphi _0\| _{\ell ^2}=1.$ Similarly $\dim \ker
(H-E_1)= 1 $ in $\ell ^2$. We denote by $\varphi _1(n)$ a generator
of $\ker (H-E_1) $ normalized so that $\| \varphi _1\| _{\ell
^2}=1.$ We   pick $\varphi _0(n)>0$ and $\varphi _1(n)\in \Bbb R$
for all $n$. Consider now the
 discrete nonlinear Schr\"odinger equation (DNLS)

$$i\partial _t u (t,n) -  (Hu )(t,n) +     |u(t,n)|^ 6u(t,n)=0
    .\tag 1.2$$
We consider a family of ground states  solutions $e^{ i\omega t}\phi
_\omega$ of (1.2), or equivalently of
$$    (Hu )( n) -     |u( n)|^ 6u( n)=-\omega  u ( n)
.\tag 1.3$$ Specifically  we have, see  Appendix A \cite{CT}:

\proclaim{Lemma 1.1} Assume (H1)--(H3). There is a family $\omega
\to \phi _\omega $ of standing waves   solving (1.3) with the
following properties. For any $\sigma \ge 0$ there is an $\eta
>0$ such that $\omega \to \phi _\omega  $ belongs to
  $C^\omega (] E_0, E_0+\eta [, \ell ^{2,\sigma})
  \cap C^0 ([E_0, E_0+\eta  [, \ell ^{2,\sigma})$.
  We have $\phi _\omega (n)
\in \Bbb R$ for any $n$ and  there are   fixed $a>0$ and $C
>0$ such that $|\phi _\omega (n)|\le C e^{-a|n|}$.
As $\omega  \searrow E_0$   we have in $C^\infty  (] E_0, E_0+\eta
[, \ell ^{2,\sigma})
  \cap C^0 ([E_0, E_0+\eta  [, \ell ^{2,\sigma})$
   the expansion
$$\phi _\omega =(\omega - E _0)^{\frac 1 {6}}  \|
\varphi _0  \| _{\ell ^8} ^{-\frac {4}{3}}(\varphi _0+ O(\omega - E
_0)) .$$
\endproclaim
Under our hypotheses, the ground states $e^{ i\omega t}\phi _\omega$
are by well known arguments orbitally stable in $\ell ^2$. A natural
question is whether they are also asymptotically stable. Here we
will define asymptotic stability as follows:

\proclaim{Definition 1.2} We say that an orbitally stable ground
state $e^{ i\omega _0 t}\phi _{\omega _0}$ of (1.2) is asymptoticaly
stable   if there are a $\sigma >0$ and an $\epsilon _0>0$ such that
for any  $u_0\in \ell ^{1,\sigma}$ with $ \|u_0- \phi_ {\omega _0}
\|_{\ell ^{1,\sigma} }  \le \epsilon _0 $  there are a $\omega _+$
and a real valued function $\theta (t)$ with
$$ \lim _{t\to +\infty }\| u(t)- e^{i\theta (t)}\phi_ {\omega
_+} \| _{\ell ^{\infty ,-\sigma} }=0.$$

\endproclaim
Our stating point is the following theorem, proved in \cite{CT}, for
  a weaker result see also \cite{KPS}:

\proclaim{Theorem 1.3} Suppose that $H$ satisfies (H3),  has just
one single eigenvalue and  $q\in \ell ^{1,1}$.  Then for any $\omega
_0\in ] E_0, E_0+\eta [$ there exist an $\epsilon_0>0$ and a $C>0$
such that if we pick $u_0\in \ell ^2$ with $  \|u_0- \phi_ {\omega
_0} \|_{\ell ^2}<\epsilon <\epsilon _0,$ then there exist
$\omega_+\in ( E_0, E_0 +\eta _0 ) $, $\theta\in C ^1 ( \Bbb R)$ and
$ u_+\in
 \ell ^2 $    with $|\omega _+-\omega _0|+\| u_+\| _{\ell ^2} \le C\epsilon  $ such that
 if $u(t,n)$ is the corresponding solution of (1.2) with
 $u(0,n)=u_0(n)$, then
$$\aligned &
\lim_{t\to\infty}\|u(t )-e^{i\theta(t)}\phi_{\omega_{+}} -e^{
it\Delta }u_+\|_{\ell ^2}=0 .
 \endaligned $$
 \endproclaim
Theorem 1.3 implies asymptotic stability of the standing waves $e^{
i\omega   t}\phi _{\omega  }$ in the sense of Definition 1.2.
Theorem 1.3 holds because $H$ is rather special. In fact,
consistently  with known results in the literature \cite{BV,JA} we
show what follows:

\proclaim{Theorem 1.4} Consider in Lemma 1.1 $\sigma >0$ large and
$\eta >0$ small.  Assume (H1)--(H3). {\item {(1)}} Any $\omega _0\in
] E_0, E_0+\eta [$   is not asymptotically stable. More precisely
for any $\sigma >0$ there is a sequence with $\| u _n(0)- \phi_
{\omega _0} \| _{\ell ^{2 , \sigma} }\to 0$ such that for any $n$
$$\inf _{\omega ,\gamma ,t } \| u_n(t)- e^{i\gamma }\phi_ {\omega  }
\| _{\ell
 ^{2 , -\sigma } } >0.$$
{\item {(2)}} For $\epsilon  >0$  small enough there is a fixed
$A_0(\omega _0)$ such that for $\| u_0-\phi _{\omega _0}\| _{\ell
^2} \le \epsilon
 $ we can write for all $t\in \Bbb R$ and for a $\theta \in C^1(\Bbb
 R, \Bbb R)$
$$ \align   & u(t,n) = e^{i\theta (t)} (\phi _{\omega (t)} (n)+
r(t,n)) , \, \,  \theta (t)=\int _0^t \omega (t') dt'+\gamma (t),
\tag 1.4\\& \text{ with $|\omega _0-\omega (t)| +\| r(t)\| _{\ell
^2} \le A_0(\omega _0)\epsilon $} \tag 1.5\\& \text{and $\langle \Re
r(t), \phi _{\omega (t)}\rangle =\langle \Im r(t),
\partial _\omega \phi _{\omega (t)}\rangle =0 $.} \tag
1.6\endalign$$ {\item {(3)}} We have a  representation  for all
times $t$
$$u(t)= e^{i\theta (t)}\left (\phi _{\omega (t)} +z\xi _1(\omega
(t))+ \overline{z}\xi _2(\omega (t))\right )+ A(\omega (t), z(t))+
h(t) \tag 1.7$$ with for a fixed $C=C(\omega _0)$

$$\align & \|A(\omega (t), z(t))\| _{\ell ^{2,2}}\le
C |z(t)|^2 \tag 1.8 \\&  \lim _{t\to \infty }\| h(t) -e^{i\Delta t}
h _{+}\| _{\ell ^{2 }}=0 \text{ for a $h_+\in \ell ^2$ with $\| h_+
\| _{\ell ^2}  \le C \epsilon $ }  \tag 1.9 \\&   \| h    \|
_{L^2_t\ell ^{2,-2 }}  \le C \epsilon    \tag 1.10 \\&  \| \xi
_{1}(\omega
 )-\varphi _1\|_{\ell ^1} +\| \xi _{2}(\omega  )\|_{\ell ^1} <C
\, \, (\omega -E_0) .\tag 1.11
\endalign$$
  There is another variable  $z=\zeta +\alpha
(\omega ,\zeta ) +\Cal A (\omega ,\zeta , h )  $ with $|\alpha
(\omega ,\zeta )|\le C| \zeta |^2$ and $|\Cal A (\omega ,\zeta , h
)|\le C | \zeta | \| h(t)\| _{\ell ^{2,-2}}$  such that
$$i\dot \zeta  -\lambda (\omega )\zeta =d(\omega ,|\zeta |^2 )\zeta
 + \Cal B (\omega ,\zeta ,h)  \tag 1.12$$
with $d(\omega ,|\zeta |^2 )$ real valued and $ |\Cal B (\omega
,\zeta ,h) |\le C \| h\| _{\ell ^{2,-2}}^2. $ There is a change of
variables $\omega =\varpi +\alpha (\varpi  ,\zeta ) +\Cal C (\varpi
,\zeta ,h)  $ with $|\alpha (\varpi ,\zeta )|\le C| \zeta |^2$ and
$\|  C (\varpi ,\zeta ,h) \| _{\ell ^{2, 2}} \le C | \zeta |  \| h
\| _{\ell ^{2,-2}}$  such that
$$i\dot \varpi    =
D(\varpi ,\zeta ,h
 )    \tag 1.13$$
with   $|  D(\varpi ,\zeta ,h
 )  |\le C \| h\| _{\ell ^{2,-2}}^2$. In particular there is a $\varpi
 _+$  with $ |  \omega _0- \varpi
 _+|<C \epsilon ^2$ such that
 $$\lim _{t\to \infty} (\zeta (t), \varpi (t)) = (0,\varpi
 _+ ).\tag 1.14$$
\endproclaim
For Stricharz estimates for $h(t)$ see \S 5.  Claim (1) is not
surprising   in view of earlier work on quasiperiodc breathers of
the DNLS in the anticontinuous limit, \cite{BV,JA}. Claim (2) is
just the standard statement of orbital stability.  What to our
knowledge is new is the detailed analysis of the long time dynamics
near the ground states. We do not push our analysis enough to prove
the existence of a manifold of quasiperiodic breathers and to prove
its asymptotic stability. Yet
  maybe Claim (3) is suggestive of such a situation.

  Our viewpoint  is the same used for continuous NLS in $\Bbb R^d$
  already used in    \cite{CT}, see references therein.
The result is in some sense opposite, but for reasons consistent
with the results for continuous NLS or in \cite{CT}. The difference
between Theorems 1.3 and 1.4 is due to   differences between the
discrete spectra $\sigma _d(H)$. The case of Theorem 1.3 is simpler
because $\sigma _d(H)$ consists of just one eigenvalue and as a
consequence $\sigma _d(\Cal H _\omega )=\{  0 \}$ with $\Cal H
_\omega$ the linearization in (2.3). However, in analogy to the
continuous case, one could prove Theorem 1.3 even in cases when
$\sigma _d(H)$ is formed by an arbitrary number of eigenvalues  and
so $\sigma _d(\Cal H _\omega )\supsetneqq \{ 0 \}$. To explain why
for $H$ satisfying (H3) Theorem 1.3 fails, we need to discuss the
nonlinear coupling between discrete and continuous modes. (1.2) can
be rewritten decomposing canonically $u(t  )$ into the spectral
components associated to the spectral decomposition of
  $\Cal H_{\omega (t)}$. In
particular $u(t  )$ is expressed as a ground state plus a reminder
$R$, see \S 2. Asymptotic stability  of ground states corresponds to
the fact that $R(t)$ disperses. $R$ is decomposed in various
components associated to the spectrum of $ \Cal H_{\omega }$. Thanks
to the dispersion theory in \cite{SK,KKK,PS,CT} the continuous
spectrum can be thought as stable spectrum. The discrete spectrum
$\sigma _d(\Cal H_{\omega })$ corresponds to central directions. To
see whether or not there is asymptotic stability, we need to study
the behavior of the discrete components of $R$. Notice that in the
absence of nonlinear coupling, the discrete components of $R$ would
describe  periodic motions. In the case of the NLS in $\Bbb R ^d$,
 for any $\lambda (\omega )
 \in \sigma _d(\Cal H_{\omega
  })\backslash \{ 0 \} $ there is a fixed integer $n$ such
 that
 $n\lambda (\omega )
 \in \sigma _e(\Cal H_{\omega
  })$, with $\sigma _e$ the     continuous
 spectrum. In other words,  there
 is resonance between continuous spectrum and nonzero
 eigenvalues of $\Cal H_{\omega
  }$. This fact seems to be responsible for the asymptotic
  stability for  NLS in $\Bbb R ^d$, and of orbital instability
  in the case of standing waves with nodes which are linearly stable, see
  \cite{C}. Using bifurcation theory as in
  Lemma 1.1 it is easy to find also  examples of  DNLS in $\Bbb Z ^d$
such that there is resonance between continuous spectrum and nonzero
 eigenvalues of $\Cal H_{\omega
  }$. In these cases one expects that energy leaks slowly from the
  discrete modes to the continuous modes because of resonances, and
 that  continuous modes disperse.  In this paper hypothesis
  (H3)   assures that $\sigma _d(\Cal H_{\omega
  })\backslash \{  0\} =\{ \pm \lambda (\omega )\} $ and that
$n\lambda (\omega )
 \not \in \sigma _e(\Cal H_{\omega
  })$ for all integers $n$. This can be easily manifactured because
$\sigma (-\Delta  ) $ is a bounded set. By the absence of
resonances,  we show that discrete and continuous
  modes  decouple. The discrete modes persist  in their
  oscillatory motion. The continuous modes consist  of
  a small component confined in a bounded region of space
   and
  of a part which scatters as solutions of the linear constant
  coefficients Schr\"odinger equation.  To prove this decoupling
  we need to perform a Siegel normal forms argument.   It is likely that the techniques
  in \cite{BV} can be used to prove the existence of a manifold of
  quasiperiodic solutions, and that our argoments show that this
  manifold is a local attractor.

\bigskip
We end with some notation. Given an operator $A$ we set
$R_A(z)=(A-z)^{-1}$ its resolvent.   We will denote by $\Cal S(\Bbb
Z)$  the set of functions $f(n)$
 rapidly decreasing as $|n|\nearrow \infty$. We will denote by
 $\Cal S(\Bbb R\times \Bbb Z)$  the set of functions $f(t, n)$
 rapidly decreasing as $(t,n) $ diverges along with all the
 derivatives $\partial _t^af(t, n)$ for $a\in \Bbb N$.
We will denote by $C_0(\Bbb Z)$  the set of functions $f(n)$
 such that $f(n)=0$ for $n$ near $\infty$.
 Given two Banach spaces $X$ and $Y$, $B(X,Y)$  will be the space of
 bounded linear operators defined in $X$ with values in $Y$. Set $\text{diag}(a,b)$ for the diagonal $2\times2$ matrix with $(a,b)$ on the
diagonal. Given a matrix $A$ we say that it is real, if it has real
valued entries. We denote by $^tA$ its transpose.

\head \S 2 Linearization, modulation and set up \endhead

In Appendix B \cite{CT} it is proved:

 \proclaim{Lemma 2.1(Global well posedness)} The DNLS (1.2)
 is globally well posed, in the sense that any initial value problem
 $u(0,n)=u_0(n)$ with $u_0\in \ell ^2$ admits exactly one solution
 $u(t)\in C^\infty (\Bbb R, \ell ^2)$. The correspondence $u_0\to u(t)$
 defines a continuous map $\ell ^2\to C^\infty ([T_1,T_2], \ell ^2)$
 for any bounded interval $[T_1,T_2] $.

\endproclaim
By the  implicit function theorem   we get:

\proclaim{Lemma 2.2 (Coordinates near standing waves)} Fix $\omega
_0$ close to $E_0$. Then there are  an $\epsilon _0>0$ and a $C_0>0$
such that any $u_0\in \ell ^2$ with $\| u_0-\phi _{\omega _0}\|
_{\ell ^2} \le \epsilon \le \epsilon _0$ can be written in a unique
way in the form $  u_0=e^{i\gamma (0)} (\phi _{\omega (0)} +r(0))$
with  $|\omega _0-\omega (0)|+|\gamma (0)|+\| r(0)\| _{\ell ^2} \le
C_0\epsilon $ and with $\langle \Re r(0), \phi _{\omega (0)}\rangle
=\langle \Im r(0), \partial _\omega \phi _{\omega (0)}\rangle =0.$
The correspondence $u_0\to (\gamma (0),\omega (0),r(0))$ is a smooth
diffeomorphism.

\endproclaim
By standard arguments one has: \proclaim{Lemma 2.3 (Orbital
Stability)} For $\epsilon _0>0$ small enough in Lemma 2.2 there is a
fixed $A_0(\omega  _0)$ such that for $\| u_0-\phi _{\omega _0}\|
_{\ell ^2} \le \epsilon \le \epsilon _0$ we can  express $u(t)$ as
(1.4)--(1.6)  for all $t\in \Bbb R$.
      \endproclaim

Consider the initial datum $u_0(n)$. We consider for all $t$ the
decomposition (2.1).   When  we plough the ansatz  (2.1) in (1.2) we
obtain

$$\aligned &
  i \partial _t r (t,n)  =
 (H r ) (t,n)  +\omega (t) r(t,n) -4
  \phi _{\omega (t)} ^6(n) r(t,n)   -3
 \phi _{\omega (t)} ^6(n)  \overline{  r }(t,n) \\& + \dot \gamma (t) \phi
_{\omega (t)}(n) -i\dot \omega (t)
\partial _\omega \phi   _{\omega (t)}(n)
+ \dot \gamma (t) r(t,n)
 +
  N (r(t,n)     ) \endaligned \tag 2.1
$$
for $N (r (t,n)    )=O(r^2(t,n))$.  We set $^tR= (r,\bar r) $,
$^t\Phi = ( \phi _{\omega } , \phi _{\omega } ) $
 and we rewrite the above equation as
$$i  R _t =\Cal H _{\omega}   R +\sigma _3 \dot \gamma   R
+\sigma _3 \dot \gamma \Phi - i \dot \omega \partial _\omega \Phi
+\Cal N(R ) \tag 2.2
$$
where for $\sigma _1=\left[ \matrix  0 &
1  \\
1 & 0
 \endmatrix \right],$ $
\sigma _2=\left[ \matrix  0 &
i  \\
-i & 0
 \endmatrix \right] $ and $
\sigma _3=\left[ \matrix  1 &
0  \\
0 & -1
 \endmatrix \right]$
we have for $\Cal H_0(\omega )=\sigma _3 \left [ H  + \omega \right
]$
$$\aligned  \Cal  H_\omega =\Cal H_0(\omega )
+V_\omega \text{ with } (V_\omega R) (n)=
 \left [-  4\sigma _3
        +3i \sigma _2
   \right ] \phi _{\omega  } ^6(n)R(n)   .\endaligned \tag 2.3$$
We have the following spectral information on $\Cal H_\omega $ :
\proclaim{Lemma 2.4} Pick $\eta >0$ in Lemma 1.1 very small. Then
for any $E_0<\omega < E_0+\eta $ we have the following facts:

{\item {(1)}} For $J$ the operator defined by $J
{^t(u_1,u_2)}:={^t(\overline{u}_1,\overline{u}_2)}$ we have
$$\sigma _1\Cal H _\omega =-\Cal H _\omega \sigma _1 , \,
[J, \Cal H _\omega ]=0, \, \sigma _3\Cal H _\omega = \Cal H _\omega
^* \sigma _3.
$$
{\item {(2)}}  The spectrum $\sigma (\Cal H _\omega )$ is symmetric
with respect to the coordinate axes.

{\item {(3)}}  The essential spectrum is $\sigma _e(\Cal H_\omega )
=\sigma _e(\Cal H_0(\omega ) ) =(-4-\omega ,-\omega ) \cup (\omega ,
4+\omega ). $

{\item {(4)}} For the generalized kernel we have $N_g(\Cal H_\omega
) = \text{span}\{  \sigma _3 \Phi _\omega  ,  \partial _\omega \Phi
_\omega \} .$

{\item {(5)}}  For the discrete spectrum  we have $\sigma _d(\Cal
H_\omega )=\{ 0,\lambda _1(\omega ), -\lambda _1(\omega ) \} $ with,
for a fixed $C$
 independent from  $\omega$,
 $|\lambda  ( \omega )- E_1-\omega | <C\, (\omega -E_0) $.  We have $\dim
 \ker(\Cal H_\omega \pm
 \lambda  (\omega ))=1$ and $N_g (\Cal H_\omega \pm
 \lambda  (\omega ))=\ker(\Cal H_\omega \pm
 \lambda  (\omega )).$

\endproclaim
{\it Proof.} (1) follows from the definition of $\Cal H_\omega$ and
the fact that $\phi _\omega (n)\in \Bbb R$. (2) follows from (1).
(3) follows from Weil's essential spectrum theorem, Theorem XIII.14
\cite{RS}. $N_g(\Cal H_\omega ) \supseteq\text{span}\{  \sigma _3
\Phi _\omega ,
\partial _\omega \Phi _\omega \}  $ follows from computation. The
equality is a consequence of elementary perturbation theory and the
fact that $\Cal H_0(\omega )$ has simple eigenvalues in $\pm (\omega
-E_0)$ and no other eigenvalues near 0. The existence of $\lambda
(\omega )  $ satisfying the above estimates follows from
perturbation theory and the fact that    $\Cal H _0(\omega )$ has
two simple eigenvalues $\pm (E_1+\omega )$. The fact that $\dim
\ker(\Cal H_\omega \pm
 \lambda  (\omega ))=1$ and $N_g (\Cal H_\omega \pm
 \lambda  (\omega ))=\ker(\Cal H_\omega \pm
 \lambda  (\omega ))$ follow from  elementary perturbation theory.
Finally, for $\eta $ small enough, the hypotheses (H2) and (H3)
imply that $\Cal H_\omega$ does not have other eigenvalues. This
ends the proof of Lemma 2.4.

\bigskip
We set $\ell ^2_d(\Cal H _{\omega} )=   N_g(\Cal H _{\omega}   )
\oplus \oplus _{ \pm  }N(\Cal H _{\omega}  \mp \lambda (\omega )) $
and $\ell _c^2(\Cal H _{\omega} )=\left \{\ell _c^2(\Cal H
_{\omega}^* )\right \} ^{\perp}.$ We have  the
 $\Cal H _{\omega}$   invariant
Jordan block decomposition
$$  \ell ^2= N_g(\Cal H _{\omega} ) \oplus \big (
\oplus _{ \pm  }N(  \Cal H _{\omega}\mp \lambda  (\omega )) \big )
\oplus \ell _c^2(\Cal H _{\omega})=N_g(\Cal H _{\omega})\oplus
 N_g^\perp (\Cal H _{\omega}^\ast  ) .\tag 2.4
$$
We have:

\proclaim{Lemma 2.5} $V_\omega \in \ell ^{p,\sigma}$, for any $ p\in
[1,\infty ]$ and $\sigma \in \Bbb R$, and $\lambda (\omega ) $
depend analytically on $\omega \in ]E_0,E_0+\eta [$. It is possible
to choose generators $\xi (\omega ) \in \ker (\Cal H _{\omega}
-\lambda (\omega ))$ so that $\xi (\omega ) \in \ell ^{p,\sigma}$ ,
for any $ p\in [1,\infty ]$ and $\sigma \in \Bbb R$,  depend
analytically on $\omega \in ]E_0,E_0+\eta [$ and we have the
normalization $\langle \xi (\omega ),\sigma _3 \xi (\omega )\rangle
=1.$ Additionally we can choose $^t\xi (\omega )= (\xi _{1}(\omega
),\xi _{2}(\omega ))$ so that for any $\tau \in \Bbb R$ there is
$C_\tau $ such that $\| \xi _{1}-\varphi _1\|_{\ell ^{1,\tau }}
<C_\tau (\omega -E_0)$ and $\| \xi _{2}\|_{\ell ^{1,\tau }} <C_\tau
(\omega -E_0)$, where $\varphi _1$ is a normalized generator of
$\ker (H-E_1)$.
\endproclaim
{\it Proof.} The fact that $V_\omega \in \ell ^{p,\sigma}$ depends
analytically on $\omega$  for $|\omega -\omega _0|\le \alpha _0$ for
some $\alpha _0>0$ follows by formula (2.3) and Lemma 1.1. Then
$\omega \to \Cal H _\omega $ is analytic from $|\omega -\omega
_0|\le \alpha _0$ into $B (\ell ^{p,\sigma},\ell ^{p,\sigma})$ for
all $p\ge 1$ and $\sigma \in \Bbb R$. In particular, $\Cal H_\omega
$ is an analytic family of operators in the sense of Kato, p. 14
\cite{RS} vol.IV (that is for any fixed $\omega _0$ and $z\not \in
\sigma (\Cal H_{\omega _0} )$, then $R_{\Cal H_\omega}(z)$ is
analytic in $\omega$  for $\omega$ sufficiently close to $\omega
_0$). Then by the Kato-Rellich theorem,  Theorem XII.8 \cite{RS},
  $\lambda (\omega ) $ depends analytically on
$\omega$. Furthermore, we have  the projection operator
$$P_{\ker (\Cal  H_\omega -\lambda (\omega ))}=\frac{i}{2\pi }
\int _{|z-E_1-\omega _0|= a_0 } R_{\Cal H_\omega}(z) dz   ,\, a_0
=\frac{E_1-4}{2}.
$$
   So $P_{\ker (\Cal H_\omega -\lambda (\omega ))}$ depends
analytically on $\omega $.
$$\aligned & P_{\ker (\Cal  H_\omega -\lambda (\omega ))}=
P_{\ker (\Cal  H_0 -E_1-\omega )}+T(\omega ) =\\& =
\frac{i}{2\pi } \int _{|z-E_1-\omega _0|= a_0 } \left (R_{\Cal H_0
}(z) -R_{\Cal H_0 }(z)V_\omega   (1+ R_{\Cal H_0 }(z)V_\omega )^{-1}
R_{\Cal H_0 }(z)\right ) dz  ,\endaligned
$$    with, for any $p,q\in [1,\infty ]$ and
$\tau , \sigma \in \Bbb R$, $\| T(\omega )\| _{B(\ell ^{p,\tau },
\ell ^{p,\sigma })}\lesssim (\omega - E_0)$ with a   fixed constant.
This implies also the information for $\xi$
 and concludes Lemma 2.5.

\bigskip
\noindent The conditions $\langle \Re r(t), \phi _{\omega
(t)}\rangle =\langle \Im r(t), \partial _\omega \phi _{\omega
(t)}\rangle =0 $ are the same of $\langle R(t), \Phi _{\omega
(t)}\rangle =\langle R(t), \sigma _3\partial _\omega \Phi _{\omega
(t)}\rangle =0 $, that is $R(t)\in N_g^\perp (\Cal H^\ast _{\omega
(t)}) $. This in particular implies that, in correspondence to the
spectral decomposition (2.5) below, we have a decomposition
$$    R (t) =(z  \xi + \bar z  \sigma _1 \xi ) + f(t)   \in
\big [ \sum _{ \pm  } \ker (\Cal H _{\omega (t)}\mp \lambda  (\omega
(t)))\big ] \oplus L_c^2(\Cal H _{\omega (t)})  .\tag 2.5  $$
$R(t)\in N_g^\perp (\Cal H^\ast _{\omega (t)}) $ implies, for $P
_{N_g(\Cal H_\omega )}$ the projection on   $N_g(\Cal H_\omega )$
associated to (2.4),
$$\aligned & i\dot \omega \langle \Phi  , \partial
_\omega \Phi
 \rangle
=\langle \sigma _3 \dot \gamma    R+\Cal N(R)+i\dot \omega \partial
_\omega P _{N_g(\Cal H_\omega )}R, \Phi \rangle
\\& \dot \gamma \langle \Phi  , \partial _\omega \Phi
 \rangle
=-\langle \text{same as above}, \sigma _3 \partial _\omega \Phi
\rangle .
\endaligned  \tag 2.6$$
We have: \proclaim{Lemma 2.6} There are two functions  $\mu (\omega
, r,\overline{r})$ and   $\nu (\omega , r,\overline{r})$ defined in
the subset of $\ell^2\times \ell ^2$ defined by $|z|+|\overline{z}|+
\| f\| _{\ell ^{2,-\tau}}\le  ( \omega -E_0)^{\frac{1}{6}} \alpha
_0$ and by $|\omega -\omega _0|\le \beta _0$ for some fixed small
$\alpha _0>0$ and $0<\beta _0\ll (\omega_0-E_0)$, analytic in
$(\omega , r,\overline{r})$, such that for any $\tau \ge 0$ we have
  $|\nu ( \omega ,r,\overline{r})|+|\mu (\omega ,r,\overline{r})|\le
  C(\tau ,\omega )  \| (r,\overline{r})\| _{\ell ^{2,-\tau}}^2$ with
  $C(\tau ,\omega ) $ continuous and with for real $\omega$
 $$\aligned & i\dot \omega =i\dot \omega (\omega ,r,\overline{r}) =
 \nu (\omega ,r,\overline{r}) -\nu (\omega ,\overline{r}, {r})
  \text{ with $\overline{\nu (\omega ,r,\overline{r})}=\nu (\overline{r}, {r})$} \\&
 \dot \gamma =\dot \gamma (\omega ,r,\overline{r})=
 \mu (\omega ,r,\overline{r}) +\mu (\omega ,\overline{r}, {r})
 \text{ with $\overline{\mu (\omega ,r,\overline{r})}=\mu (\omega ,\overline{r},
 {r})$}.\endaligned \tag 2.7
 $$

\endproclaim
{\it Proof.} Apply $ P_{N_g(\Cal H_\omega )}$  to (2.2) obtaining
$$i  P_{N_g(\Cal  H_\omega )}R _t - \dot \gamma (t) P_{N_g(\Cal  H_\omega )}\sigma _3   R
-\sigma _3 \dot \gamma (t)\Phi _\omega+i \dot \omega (t)\partial
_\omega \Phi _\omega+P_{N_g(\Cal  H_\omega )}\Cal N(R  )   .$$ Set
$q(\omega)=\| \phi _\omega \| ^2_2$ and $q'(\omega)
=dq(\omega)/d\omega$. Then we have
$$ P_{N_g(\Cal H_\omega )} =\sigma _3\Phi _\omega \langle \quad , \sigma
_3\partial _\omega \Phi _\omega \rangle /q'(\omega )+\partial
_\omega \Phi _\omega \langle \quad ,   \Phi _\omega \rangle
/q'(\omega ) .$$ By $P_{N_g(\Cal H_\omega )}R=0 $, which implies $
P_{N_g(\Cal H_\omega )}R _t=-\dot \omega  \partial _\omega P _{ N_g
(\Cal H _\omega )} R$, we get

$$\aligned & \left (  q'(\omega) +\left [    \matrix -\langle  \partial _\omega P _{ N_g
(\Cal H   _\omega )}  R , \Phi _\omega \rangle
&    \langle \sigma _3 R , \Phi _\omega \rangle  \\
- \langle  \partial _\omega P _{ N_g  (\Cal H   _\omega )}  R,
\sigma _3
\partial _\omega \Phi _\omega \rangle  & \langle   R , \partial
_\omega \Phi _\omega \rangle
\endmatrix \right ] \right )     \left [ \matrix  i\dot \omega \\ -\dot
\gamma
\endmatrix  \right ]  = \left [ \matrix  \langle \Cal N(R  ), \Phi _\omega \rangle  \\
\langle \Cal N(R  ), \sigma _3 \partial _\omega \Phi _\omega \rangle
\endmatrix \right ] .
\endaligned  $$

By an elementary computation we have
$$ \aligned & \langle  \partial _\omega P _{ N_g
(\Cal H   _\omega )}  R , \Phi _\omega \rangle =\left \langle
r+\overline{r},  \partial _\omega {\phi _\omega }  \right \rangle
\\&\langle
\partial _\omega P _{ N_g (\Cal H   _\omega )}  R , \sigma _3
\partial _\omega \Phi _\omega \rangle =\left \langle r-\overline{r},
 {\partial _\omega ^2\phi _\omega } \right  \rangle
\endaligned \tag 2.8$$
and so $^t( i\dot \omega , -\dot \gamma )=$

$$\aligned & \left (   q'(\omega )   +\left [    \matrix
-\langle   r+\overline{r},\partial _\omega  \phi _\omega  \rangle
&   \langle r-\overline{r} , \phi _\omega \rangle  \\
- \langle r-\overline{r}, \partial _\omega ^2\phi _\omega  \rangle
& \langle r+\overline{r} , \partial _\omega \phi _\omega \rangle
\endmatrix \right ] \right ) ^{-1}       \left [ \matrix   \langle N(r,\overline{r}) -N( \overline{r},r), \phi _\omega \rangle  \\
  \langle N(r,\overline{r})+N( \overline{r},r),   \partial _\omega \phi _\omega
\rangle
\endmatrix \right ] .
\endaligned  \tag 2.9$$
Notice that the right factor in (2.9) is a polynomial in
$(r,\overline{r})$ while the left factor is analytic in the
functions in (2.8) and   $\langle r-\overline{r}, \phi _\omega
\rangle$. This completes Lemma 2.6.

\medskip

We plug  now decomposition  (2.6) in system (2.3)  and  we obtain
$$\aligned   i\dot z    -\lambda  (\omega ) z   &=
\langle
   \dot \gamma   \sigma _3R+\Cal N(R)-\\&
 -iz \dot \omega
 \partial _\omega \xi  -i\overline{z} \dot \omega
\sigma _1 \partial _\omega \xi  +i\dot \omega   \partial _\omega
 P_{\ker (\Cal H  _\omega -\lambda  )}R, \sigma _3 \xi \rangle
  \\  i  \dot f -\Cal H_\omega f &= \dot \gamma P_c(\Cal H_\omega )
    \sigma _3R+P_c(\Cal H_\omega ) \Cal
N(R)+i\dot \omega   \partial _\omega P_{c} (\Cal H _\omega )R
.\endaligned \tag 2.10$$

Our first step will consist in \S 4 in splitting $f=\Phi (\omega ,
z) +g$ with $g\in \ell ^2_c(\Cal H_\omega )$ satisfying an equation
of the form (4.1) below and with $\|\Phi (\omega , z)\| _{\ell
^{2,2}}\le C  |z|^2$ for fixed $C>0$. In \S 5 we will prove that $g$
is asymptotically free. In \S 6 we will prove that $z(t)$ does not
decay to 0. We will also show in \S 7 that $\omega (t)$ oscillates.
We first state some linear dispersive estimates needed later. Proofs
of Lemmas 3.1 and 3.3 are in sections 8--10.

 \head \S 3 Spacetime estimates for $\Cal H _\omega  $ \endhead
We list a number of linear estimates needed  later. The constants
$C(\omega)$ and $C(\tau ,\omega)$ in this section are upper
semicontinuous  in $\omega$.

\proclaim{Lemma  3.1} Under hypotheses (H1-3)   there is a constant
$C(\omega)$    such that

$$\| P_c(\Cal H_\omega )e^{it\Cal H_\omega } \|
_{B(\ell ^p ,\ell ^{p'})}
   \le C(\omega) \langle t \rangle ^{-\frac 2 3
\left (\frac 1 p-\frac 1 2\right )}  \text{   $\forall \, p\in
[1,2]$ and for $p'=\frac {p} {p-1}$.}$$
\endproclaim

The proof is in \S 10. The next estimates needed are Stricharz
estimates. Following \cite{CV} for every $1\leq p, q\leq \infty$ we
introduce the Birman-Solomjak spaces
$$
\ell ^p({   \Bbb Z}, L^q_t[n,n+1])\equiv \left \{f\in
L^q_{loc}({\Bbb R}) \hbox{ s.t. } \{\|f\|_{L^q[n, n+1]} \}_{n\in
{\Bbb Z}} \in \ell ^p({\Bbb Z})\right \}, $$  endowed with the
  norms
$$\aligned &    \| f \|_{\ell ^p({  \Bbb Z}, L^q_t[n,n+1])} ^p\equiv
\sum_{n\in \Bbb Z}   \| f \| _{L^q_t[n,n+1] }^{p}   \quad  \forall
\quad  1\le p<\infty \text{ and } 1\leq q \leq \infty  \\ & \| f \|
_{\ell ^\infty ({ \Bbb Z}, L^q_t[n,n+1])} \equiv \sup _{n\in \Bbb Z}
\| f \| _{L^q[n,n+1] }.
\endaligned $$
We  say that a pair of numbers $(r,p) $ is admissible if
 $$ 2/{r }+ 1/{p }=1/2 \text{ and }
 (r,p) \in [4, \infty]\times [2, \infty]  .\tag 3.1$$
Then  proceeding as in \cite{CV}, by a standard $TT^\ast$ argument
it is possible to prove from Lemma 3.1 the following result:

\proclaim{Lemma 3.2 (Strichartz estimates)} Under the hypotheses
of Theorem 1.3 there exists a constant $C =C(\omega ) $
 upper semicontinuous in $\omega$  such that for every admissible pair $(r,p) $   we have:
$$
  \|e^{it\Cal H_\omega  }P_c(\Cal H_\omega   ) f \|
    _{\ell ^{\frac 32 r}( \Bbb {Z},L^{\infty}_{t}([n,n+1],
    \ell ^{ p}(\Bbb {Z})))} \leq
C (\omega )  \|f\|_{\ell ^2(\Bbb {Z})}.
$$
Moreover,  for any two admissible pairs $(r_1,p_1), (r_2,p_2) $   we
have the   estimate
$$ \aligned & \left \| \int _{0}^{t}e^{i(t-s)
\Cal H_\omega  }P_c(\Cal H_\omega   ) g(s )ds\right \| _{\ell
^{\frac{3}{2}r_1}( \Bbb {Z},L^{\infty }_{t}([n,n+1],\ell ^{
p_1}({\Bbb Z}))} \le \\& \le C(\omega ) \| g \| _{\ell ^{\left (
\frac{3}{2}r_2\right ) '}( \Bbb {Z},L^{1}_{t}([n,n+1],\ell  ^{
p_2'}({\Bbb Z}))}.\endaligned
$$\endproclaim
In \S 9 we prove the following Kato smoothness result:

\proclaim{Lemma 3.3} For $\tau >1$ there exists $C=C(\tau , \omega
)$ such that  for all $z\in \Bbb C \backslash \sigma _e(\Cal
H_\omega )$
$$\| R_{\Cal H_\omega }(z )P_c(\Cal H_\omega   )   \| _{B (\ell ^{2, \tau }
,\ell ^{2,-\tau }) }\le  C   .$$ The  following limits are well
defined for any $\lambda \in [0,4]$ in $C^0([0,4],B(\ell ^{2,\tau }
, \ell ^{2,-\tau }  ))$
$$\lim _{\epsilon \to 0 ^+}  R_{\Cal H_\omega }(\lambda \pm i\epsilon )
=R^{\pm}_{\Cal H_\omega  }(\lambda  )   .$$ For any $u\in \ell
^{2,\tau}\cap \ell ^2_c(\Cal H_\omega )$ we have
$$\aligned & P_c(\Cal H_\omega  )u=\frac{1}{\sqrt{2\pi}i}\int_\Bbb R
 (R_{\Cal H_\omega }^{+}(\lambda  )-R_{\Cal H_\omega
}^{-}(\lambda  ))   u d\lambda
\\& =\frac{1}{\sqrt{2\pi}i}\int_{\sigma _e(\Cal H_\omega )}
 (R_{\Cal H_\omega }^{+}(\lambda  )-R_{\Cal H_\omega
}^{-}(\lambda  ))   u d\lambda   .\endaligned $$
\endproclaim
 The first two statements are proved in
Lemma 9.1. The third statement is proved in Lemma 9.4. We list now a
number of corollaries.

\proclaim{Lemma 3.4}For $\tau >3/2$ $\exists$ $C=C (\tau ,\omega )$
s.t.:
   {\item {(a)}}
  for any $f\in \Cal S(\Bbb Z)$,
$$\align & \| e^{-it\Cal H_\omega }P_c(\Cal H_\omega  )f\|
_{L_t^2\ell  ^{2 , -\tau }  } \le
 C \|f\|_{\ell ^2};
\endalign $$
  {\item {(b)}}
  for any $g(t,n)\in
 \Cal S (\Bbb R \times \Bbb Z)$
$$ \left\|\int_\Bbb R e^{it\Cal H_\omega }
P_c(\Cal H_\omega  )g(t,\cdot)dt\right\|_{\ell ^2 } \le
C  \| g\|_{L_t^2\ell  ^{2, \tau} }.
$$\endproclaim
 The proof is the same of Lemma 3.3 \cite{CT}.

\proclaim{Lemma 3.5} For any  $\tau >1$ $\exists$ $C=C(\tau , \omega
) $ such that
$$\align &  \left\|  \int_0^t e^{-i(t-s)\Cal H_\omega }
P_c(\Cal H_\omega )g(s,\cdot)ds\right\|_{L_t^2\ell ^{2 ,-\tau }} \le
C \| g\|_{L_t^2\ell ^{2, \tau } }.  \endalign
$$
\endproclaim
 The proof is the same of Lemma 3.4 \cite{CT}.

\proclaim{Lemma 3.6}  For every $\tau >3/2$  $\exists$ $C=C(\tau
,\omega ) $ such that
$$
\left\|\int_0^t e^{-i(t-s)\Cal H_\omega }P_c(\Cal H_\omega
)g(s,\cdot)ds \right\|_{ L_t^\infty \ell  ^2\cap \ell ^{6}( \Bbb
{Z},L^{\infty }_{t}([n,n+1],\ell ^\infty ))} \le C \|g\|_{L_t^2\ell
^{2,\tau }}.$$
\endproclaim
The proof is the same of Lemma 3.5 \cite{CT}.
 We will now assume Lemmas 3.1 and 3.3
and we will proceed with the proof of our nonlinear result.

\head \S 4  Decoupling between localized and dispersive radiation
\endhead

\proclaim{Lemma 4.1} There is a representation $ f=\Phi (\omega , z)
+g$ such that $\sigma _1g= \overline{g}$,    $\Phi (\omega , z)\in
\ell ^2_c(\Cal H _\omega )$ is analytic in $(\omega
,z,\overline{z})$ with $\|\Phi (\omega , z)\| _{\ell ^{2,2}}\le C
|z|^2$ for fixed $C>0$, and $g\in \ell ^2_c(\Cal H _\omega )$ which
satisfies, for $\| \widehat{\Cal F} (\omega, z, g )\| _{\ell ^{2,2 }
} \le C |z|\, \| g \| _{\ell ^{2,-2 } }$,

$$\aligned & i  \dot g -
(\Cal H_\omega +\dot \gamma P_c(\Cal H_\omega )  \sigma _3) g  =
  \widehat{\Cal F}
(\omega, z, g )+O(|g|^7) .\endaligned \tag 4.1$$

\endproclaim
{\it Proof.} We enter the splitting (2.6) in equations (2.7) and
(2.10). We set $f_1=f$ and setting  $\dot \omega = \dot \omega
(\omega , R )$  $\dot \gamma = \dot \gamma (\omega , R )$, i.e. the
functions in Lemma 2.6, we write
$$\aligned &
 i \dot f_1 -\Cal H_\omega f _1 =
 \dot \gamma  P_c(\Cal H_\omega )
    \sigma _3f_1+ A_1(\omega , z )+ {\Cal F}_1(\omega , z, f_1)\\& i\dot \omega =b_1(\omega , z
)+ {\Omega} _1(\omega , z,f_1)
  \\&
 \dot \gamma =c_1(\omega , z )
 + {\Gamma} _1(\omega , z,f_1)   \\&
 i\dot z    -\lambda  (\omega ) z =d_1(\omega , z )+ {Z} _1(\omega , z,f_1)
  . \endaligned  \tag 4.2$$
We will define recursively a sequence of systems

$$\aligned & i \dot f_\ell  -\Cal H_\omega f _\ell  =\dot \gamma
  P_c(\Cal H_\omega )
     \sigma _3f_\ell + A_\ell(\omega , z
)+ {\Cal F}_\ell (\omega , z, f_\ell) \\& i\dot \omega =b_\ell
(\omega , z )+ {\Omega} _\ell (\omega , z,f_\ell )   \\&
 \dot \gamma =c_\ell (\omega , z )+ {\Gamma } _\ell (\omega , z,f_\ell )
   \\&
 i\dot z    -\lambda  (\omega ) z =d_\ell (\omega , z )
 + {Z} _\ell (\omega , z,f_\ell )
   \endaligned  \tag 4.3 $$
and we will assume a number of { \bf inductive hypotheses} for fixed
$\varepsilon _0>0$ small enough and $C$:

{\item {(1)}} $f_\ell  \in \ell ^2_c(\Cal H _\omega )$;

{\item {(2)}} ${\Cal F}_\ell (\omega , z, f_\ell )= {\Cal F}  (
f_\ell ) +\widehat{{\Cal F}}_\ell (\omega , z, f_\ell)$ with ${\Cal
F}  ( f_\ell ) =O(|f_\ell |^7)$ and with $\widehat{{\Cal F}}_\ell
(\omega , z, f_\ell )\in \ell ^{2,2}$;

{\item {(3)}}  $b_\ell $,  $c_\ell $, $d_\ell $, $ \Omega _\ell $,
$\Gamma _\ell $  and $Z_\ell $ are analytic functions in $(\omega ,
z, \overline{z}, f_\ell )$ for

  $$\max \{ |\omega -\omega _0|,|z|, \| f_\ell \| _{\ell
^{2,-2}}\} \le e^{-(2-2^{-\ell})} \varepsilon _0;\tag 4.4$$
furthermore, their Taylor expansions  have real coefficients;

{\item {(4)}} $(\omega , z,f_\ell ) \to \widehat{{\Cal F}}_\ell
(\omega , z, f_\ell)\in \ell ^{2,2}$
 is an  analytic  function in the domain in (4);

{\item {(5)}} $(\omega , z  ) \to A_\ell(\omega , z )\in \ell
^{2,2}$
 is   analytic   in $(\omega
, z, \overline{z}, f_\ell )$ for $\max \{  |\omega -\omega _0|,|z|
\} \le e^{-(2-2^{-\ell})} \varepsilon _0;$ if we expand

$$A_\ell(\omega , z )=
\sum _{m+n\ge \ell +1}A_{\ell m n}(\omega ) z^m \overline{z}^n
\text{ and set }\widetilde{A}_\ell(\omega , z ):= \sum _{m+n\ge \ell
+2}A_{\ell m n}(\omega ) z^m \overline{z}^n,$$ then the  $A_{\ell m
n}(\omega )$ are real with $\sigma _1A_{\ell m n}(\omega )=- A_{\ell
n m}(\omega )$; $A_\ell  (\omega ,z) \in \ell ^2_c(\Cal H _\omega )
\cap  \ell ^{2,2}$;

{\item {(6)}} The subspace of solutions of (4.3) with $\omega$ real,
$z$ and $\overline{z}$ complex conjugate of each other and $\sigma
_1f_\ell =\overline{f}_\ell$, is invariant by system (4.3);

{\item {(7)}} the following estimates hold in (4.4)

$$\aligned & \|A_\ell(\omega , z ) \| _{\ell ^{2, 2}}\le
   C_A(\ell )  e^{ (\ell -1)(2-2^{-\ell})}
\varepsilon _0 ^{-\ell +1 }  |z|^{\ell +1 } \\&  \max \{ (|\Omega
_\ell |,|\Gamma _\ell |,|Z _\ell | , \| \widehat{\Cal F} _\ell
  \| _{\ell ^{2, 2}})(\omega, z ,f_\ell )\} + \le C  _\Omega (\ell )
 \left ( |z|+  \| f_\ell  \| _{\ell ^{2,
-2}}\right ) \| f_\ell  \| _{\ell ^{2, -2}}
\\& \max \{ |b_\ell (\omega , z )|,|c_\ell (\omega , z )|,|d_\ell (\omega ,
z )|\}\le C_b(\ell )  |z|^2.
\endaligned \tag 4.5$$

These hypotheses hold for $\ell =1$, for   $ C_A(1 )=C_\Omega (1
)=C_b(1 )=c(1)$.  We   define for $\Phi _{(\ell +1)mn}(\omega  ):=R
_{\Cal H _\omega }((m-n)\lambda (\omega )) A _{\ell mn}(\omega
   )$
$$f_{\ell  }=f_{\ell +1}+ \Phi _{\ell +1}(\omega , z) \, ,\,
\Phi _{\ell +1}(\omega , z)=\sum _{ m+n =\ell +1} \Phi _{(\ell
+1)mn}(\omega  ) z^m \overline{z}^n.$$ We have $\Phi  _{\ell +1} \in
\ell ^2_c(\Cal H _\omega )$. This implies
  $f  _{\ell +1}  \in \ell ^2_c(\Cal H _\omega )$.
We have:

\proclaim{Lemma 4.2} $ \Phi _{\ell +1}(\omega , z)$ is analytic in
$(\omega , z,\overline{z})$ for  $|\omega -\omega _0|\le
e^{-(2-2^{-\ell   })} \varepsilon _0$ and $ (z,\overline{z})\in \Bbb
C^2 $ with values in $\ell ^{2,2}$.
  In $|\omega
-\omega _0| \le e^{-(2-2^{-\ell -1})} \varepsilon _0  $ and   for $
\Phi _{\ell +1}' \cdot  (b , c ):=\partial _\omega \Phi _{\ell +1}
b+\partial _z \Phi _{\ell +1} c-\partial _{\overline{z}}
\overline{\Phi} _{\ell +1} \overline{c} $ we have for $D_\ell
=C_{4.9}\frac{ e^{2 } 2^{2\ell +2}}{\ell !} $, see (4.9) for
$C_{4.9}$,
$$\aligned  & \|\Phi _{\ell +1}(\omega , z)\| _{\ell ^{2,2}}   \le
C_A(\ell )D_\ell
 e^{ (\ell
-1)(2-2^{-\ell})} \varepsilon _0 ^{-\ell +1 } |z|^{\ell +1 }\le
C_A(\ell )D_\ell   |z|^{2
  }
\\&  \|  \Phi _{\ell +1}'(\omega , z)\| _{\ell ^{2,2}}  \le 4
C_A(\ell )D_\ell e^{ (\ell -1)(2-2^{-\ell})} \varepsilon _0 ^{-\ell
+1 } |z|^{\ell
  }\le 4C_A(\ell )D_\ell   |z|
 .\endaligned
 \tag 4.6
$$
 For $\omega $ real   $\Phi
_{ (\ell +1)mn}(\omega  )$  is real,  $\sigma _1\Phi _{(\ell
+1)mn}(\omega  )=\Phi _{(\ell +1)nm}(\omega  )$, and so
 $\sigma _1f _{\ell +1}
 =\overline{ f _{\ell +1}  }$ if $\sigma _1f _{\ell  }
 =\overline{ f _{\ell   }  }$. In particular this yields  the
 inductive hypothesis (6) for $\ell +1$.
\endproclaim
{\it Proof.} The last  two sentences  follow   from $\sigma
_1A_{\ell mn}(\omega  )= - A_{\ell mn}(\omega  )$  for $\omega $
real and from $\sigma _1\Cal H _\omega = -\Cal H _\omega \sigma _1.$
Now we turn to the estimates. There are fixed $C_0>0$ and $\alpha
_0>0$ such that for $|\omega -\omega _0|\le
 \varepsilon _0$ and $j=0,1$
$$|\partial _{\omega }^j\left [R _{\Cal H   _0(\omega )
}((m-n)\lambda (\omega ), \mu , \nu )\right ] |\le C_0 e^{-\alpha
_0|\mu -\nu|} .\tag 4.7$$ Set $P_c(\omega )=P_c(H _\omega )$ and $
\Cal H _\omega  ^c = P_c(\ \omega ) \Cal H _\omega  $. Then
$$\aligned & R_{\Cal H _\omega  ^c} =R _{\Cal H   _0(\omega )}P_c(\omega )
-R _{\Cal H   _0(\omega )}   V_\omega P_c(\omega ) R _{\Cal H
_0(\omega )} +R _{\Cal H _0(\omega )} V_\omega R_{\Cal H _\omega ^c}
V_\omega  R _{\Cal H _0(\omega )}.
\endaligned \tag 4.8$$
By (4.7)--(4.8) for $|\omega -\omega _0|\le \varepsilon _0$ and for
all $(m,n )\in \Bbb Z^2$
$$\|\partial _\omega ^j \left [ R _{\Cal H ^c_\omega
}((m-n)\lambda (\omega ) )\right ] \| _{B(\ell ^{2,2}, \ell
^{2,2})}\le C_{4.9} \text{ for $j=0,1$} .\tag 4.9$$ (4.9) implies
for $\max \{ |\omega -\omega _0|,|z|   \} \le e^{-(2-2^{-\ell-1})}
\varepsilon _0 $

$$\aligned & \|\Phi _{\ell +1} \| _{\ell ^{2,2}}\le
 C_{4.9}  \sum _{m+n=\ell +1}\frac{1}{m!n!} \|  \partial _z^m
\partial _{\overline{z}}^nA_\ell (\omega , 0)\| _{\ell ^{2,2}}  \, |z|^{\ell
+1}\\& \le       \frac{ C_{4.9}    2^{\ell +1}|z|^{\ell +1} }{(\ell
+1)! \left ( \varepsilon _0 \exp (2^{-\ell}-2 )\right )^{\ell +1}}
\sup _{ |z| \le e^{-(2-2^{-\ell })} \varepsilon _0} \| A_\ell
(\omega , z)\| _{\ell ^{2,2}} \\& \le  C_{4.9}C_A(\ell ) \frac{
2^{\ell +1}|z|^{\ell +1} }{(\ell +1)! \left ( \varepsilon _0 \exp
(2^{-\ell}-2 )\right )^{\ell -1}}     .
\endaligned \tag 4.10
$$
 Differentiating we obtain
  $\partial _\omega \Phi _{\ell +1} =
  \Phi _{\ell +1}^{(1)}+\Phi _{\ell +1}^{(2)}$
   with

$$ \aligned & \Phi _{\ell +1}^{(1)}(\omega , z )=
\sum _{ m+n =\ell +1}  \partial _\omega \left [ R _{\Cal H ^c
_\omega }((m-n)\lambda (\omega )  )\right ] A _{\ell mn}(\omega
    ) z^m \overline{z}^n\\& \Phi _{\ell +1}^{(2)}(\omega , z )=
\sum _{ m+n =\ell +1}  R _{\Cal H ^c _\omega }((m-n)\lambda (\omega
)  ) \partial _\omega A _{\ell mn}(\omega
    ) z^m \overline{z}^n.\endaligned$$
We have $\|\Phi _{\ell +1}^{(1)}(\omega , z )\| _{\ell ^{2,2}}\le
C_{4.9} C _A(\ell ) \frac{      2^{\ell +1}|z|^{\ell +1} }{(\ell
+1)! \left ( \varepsilon _0 \exp (2^{-\ell}-2 )\right )^{\ell -1}}
$ as for (4.10). By
$$  \varepsilon _0 \exp (2^{-\ell }-2 ) -
 \varepsilon _0 \exp (2^{-\ell -1}-2 )\ge \varepsilon _0 e^{-2}
2^{-\ell -1}  \tag 4.11$$   and by the Cauchy integral formula, we
have

$$\aligned & \|\Phi _{\ell +1}^{(2)}(\omega , z )\|
_{\ell ^{2,2}}\le C_{4.9}    \sum _{m+n=\ell +1}\frac{1}{m!n!}
\|\partial _z^m
\partial _{\overline{z}}^n\partial _\omega A_\ell (\omega , 0)
\| _{\ell ^{2,2}}  \, |z|^{\ell +1}\\& \le  C_{4.9}  C_A (\ell )
\frac{
  2^{ 2\ell +2} e^{ 2}|z|^{\ell +1}}{(\ell +1)! \varepsilon _0 \left (
\varepsilon _0 \exp (2^{-\ell}-2 )\right )^{\ell -1}} .
\endaligned
$$
So we conclude

$$\aligned & \|\partial _\omega \Phi _{\ell +1}(\omega , z )\|
_{\ell ^{2,2}}\le 2 C_{4.9}  C _A (\ell )\frac{
  2^{ 2\ell +2} e^{ 2}|z|^{\ell +1}}{(\ell +1)! \varepsilon _0 \left (
\varepsilon _0 \exp (2^{-\ell}-2 )\right )^{\ell -1}} .
\endaligned \tag 4.12
$$

We have $$\aligned & \| \partial _z\Phi _{\ell +1} \| _{\ell
^{2,2}}\le
 C_{4.9}   \sum _{m+n=\ell +1}\frac{1}{(m-1)!n!}  \|\partial _z^{m }
\partial _{\overline{z}}^nA_\ell (\omega , 0)\| _{\ell
^{2,2}}
  \, |z|^{\ell
 }\\& \le    C_{4.9}  C _A(\ell ) \frac{
  2^{  \ell +1}  |z|^{\ell  }}{ \ell  !  \left (
\varepsilon _0 \exp (2^{-\ell}-2 )\right )^{\ell -1}} .
\endaligned \tag 4.13
$$
     This concludes the proof of Lemma 4.2.

\bigskip

 We get equations
(4.3) for $\ell +1$
 with, for ${\Omega }_{\ell +1} ( f_{\ell +1})=
 {\Omega }_{\ell +1} (\omega, z, f_{\ell +1})$ etc.,
$$\aligned & b_{\ell +1} (\omega , z )= b_{\ell  } (\omega , z )
+
  {\Omega}  _{\ell  } (\omega ,
z,\Phi _{\ell +1} )  \, , \,   {\Omega }_{\ell +1} ( f_{\ell +1}) =
 {\Omega} _{\ell  } ( f_{\ell +1} +\Phi _{\ell +1}) -
 {\Omega} _{\ell  } ( \Phi _{\ell +1})\\& c_{\ell +1} (\omega
, z )= c_{\ell } (\omega , z ) +
  {\Gamma}  _{\ell  } (\omega ,
z,\Phi _{\ell +1} )  \, , \,  {\Gamma} _{\ell +1} ( f_{\ell +1}) =
 {\Gamma} _{\ell } ( f_{\ell +1} +\Phi _{\ell +1})-
 {\Gamma} _{\ell } ( \Phi _{\ell +1})
\\&  d_{\ell +1} (\omega , z )=
d_{\ell } (\omega , z ) + \
  {Z} _{\ell  } (\omega ,
z,\Phi _{\ell +1} )  \, , \,    {Z} _{\ell +1} ( f_{\ell +1}) =  {Z}
_{\ell  } ( f_{\ell +1} +\Phi _{\ell +1})- {Z} _{\ell } ( \Phi
_{\ell +1})  ,
\endaligned  $$

$$\aligned &A_{\ell +1}  = P_c(\Cal H _\omega ) \big
\{\widetilde{A}_\ell +c_{\ell +1} \sigma _3 \Phi _{\ell +1}    +
\Cal F _{\ell } (\omega , z, \Phi _{\ell +1}) -\Phi _{\ell +1}'
\cdot (b_{\ell+1}, c_{\ell+1} ) \big \} ,
\endaligned  \tag 4.14$$

$$\aligned &  {\Cal F}_{\ell +1}
(\omega , z, f_{\ell +1}) =  P_c(\Cal H _\omega )
 \big \{  {\Cal F}_\ell (  f_{\ell +1}
+\Phi _{\ell +1})   - {\Cal F}_\ell (  \Phi _{\ell +1}) \\&  +
\Gamma _{\ell+1}( f_{\ell +1})  \sigma _3\Phi _{\ell} -\Phi _{\ell
+1}' \cdot \left (   \Omega _{\ell +1}( f_{\ell +1}) , Z _{\ell +1}
( f_{\ell +1})\right ) \big \} +i\dot \omega
\partial _\omega  P_c(\Cal H _\omega ) f _{\ell +1} .
\endaligned \tag 4.15$$
Notice that $\sigma _1 A_{\ell +1,m n}(\omega ) =- { A_{\ell +1,
nm}}(\omega )$ for  $\omega$  real can be derived from the
invariance  of (4.3) for $\ell +1$ solutions of the space where
$\omega $ is real, $z$ and $\overline{z}$ are complex conjugates and
$\sigma _1 f_{\ell +1}= \overline{ f_{\ell +1}}$, already stated in
Lemma 4.2.

\proclaim{ Lemma 4.3} The following, are analytic functions in $(
\omega , z,\overline{z}, f_{\ell +1})$:

{\item {(1) } }$ b_{\ell +1} ,$ $ c_{\ell +1}
 $ and  $ d_{\ell +1}  $ are analytic in $\max \{ |\omega -\omega _0|
 ,|z|
  \}
\le e^{-(2-2^{-\ell -1})}\varepsilon _0$; {\item {(2) } }  $A_{\ell
+1} $ is analytic in $\max \{ |\omega -\omega _0|,|z| \}\le
e^{-(2-2^{-\ell -1})}\varepsilon _0$  with values in $\ell ^{2,2}$;
{\item {(3) } } for $\max \{ |\omega -\omega _0|,|z|,
   \| f_{\ell +1}\| _{\ell ^{2,-2}}\} \le
e^{-(2-2^{-\ell -1 })}\varepsilon _0$,  $  {\Omega }_{\ell +1}$, $
{\Gamma }_{\ell +1}$ and $  {Z }_{\ell +1}$ are    analytic; {\item
{(4) } } we have $ {\Cal F}_{\ell +1}(\omega , z, f_{\ell +1})=
{\Cal F} ( f_{\ell +1}) + \widehat{\Cal F}_{\ell +1} (\omega , z,
f_{\ell +1}) $ with ${\Cal F}$ the same of induction hypothesis (2)
and with $ \widehat{\Cal F}_{\ell +1} (\omega , z, f_{\ell +1}) $
analytic  with values in $\ell ^{2,2}$ and with domain $\max \{
|\omega -\omega _0|,|z|,
  \| f_{\ell +1}\| _{\ell ^{2,-2}}\}  \le
 e^{-(2-2^{-\ell-1})}\varepsilon _0$.
\endproclaim

{\it Proof.} (1)--(2) are   consequences of $C_A(\ell )D_\ell e^{
 (2^{-\ell  -1} -2 )} \varepsilon _0  < 1 $ which we assume,
 see (4.27). Indeed

$$ \|\Phi _{\ell +1}(\omega , z)\| _{\ell ^{2,2}}
\le C_A(\ell )D_\ell
 e^{ (\ell
-1)(2-2^{-\ell})} \varepsilon _0 ^{-\ell +1 } |z|^{\ell +1 }\le
C_A(\ell )D_\ell e^{  (2^{-\ell  } -4)} \varepsilon _0 ^{2 }  $$
implies $\max \{ |\omega -\omega _0|,|z|,
  \| \Phi _{\ell +1}(\omega , z)\| _{\ell ^{2,-2}}\}  \le
 e^{-(2-2^{-\ell-1})}\varepsilon _0$.  For (3), write $$\| \Phi _{\ell +1}\| _{\ell ^{2,-2}} +\| f _{\ell +1}\| _{\ell
^{2,-2}}\le C_A(\ell )D_\ell |z|^{2}+  e^{-(2-2^{-\ell -1} )}
\varepsilon _0\le e^{-(2-2^{-\ell  } )}   \varepsilon _0$$ The  last
inequality follows from (4.11) and the following inequality, which
we assume, see (4.27),
$$C_A(\ell )D_\ell |z|^{2}\le C_AD_\ell  \varepsilon _0^{2}
\le \varepsilon _0 e^{-2}2^{-\ell -2 }.\tag 4.16$$
  To prove (4) we proceed similarly by
$  \Cal F _{\ell +1}( \omega , z ,f _{\ell +1}) =\Cal F  (  f _{\ell
+1}) +\widehat{\Cal F} _{\ell +1}( \omega , z ,f _{\ell +1})
 $

$$\aligned &\widehat{\Cal F} _{\ell +1}( \omega , z ,f _{\ell +1})
= P_c(\Cal H _\omega )\widehat{\Cal F} _{\ell  }( \omega , z ,f
_{\ell +1 })
\\& +\int _{[0,1]^2} P_c(\Cal H _\omega )\Cal F _{\ell  }''( \omega , z ,sf _{\ell
}+t\Phi _{\ell +1} )dt\, ds\, \cdot f _{\ell +1 }\cdot \Phi _{\ell
+1}  +\\&    P_c(\Cal H _\omega )
 \big \{  \Gamma _{\ell+1}( f_{\ell +1})
\sigma _3\Phi _{\ell +1 }  -\Phi _{\ell +1}'  \cdot \left ( \Omega
_{\ell +1} ,   Z _{\ell +1} \right ) \big \} +i\dot \omega
\partial _\omega  P_c(\Cal H _\omega ) f _{\ell +1}
\endaligned \tag 4.17
$$

 \proclaim{Lemma 4.4} In (4.5) we can choose $C_b $, $C_\Omega  $
 and  $C_A $ in $\ell ^\infty (\Bbb N)$.
\endproclaim
{\it Proof.} We have $|\Omega _{\ell } (\omega , z,\Phi _{\ell +1} )
| \le 2C_\Omega (\ell ) |z| \| \Phi _{\ell +1} \| _{\ell ^{2,2}}\le
2C_\Omega (\ell )C_A(\ell )  D_\ell   |z| ^{3}   $ for $C_A(\ell
)D_\ell \varepsilon _0<1$. So for $C_{18}(\ell ):=   C_b (\ell )+
2\varepsilon _0  C_\Omega (\ell ) C_A (\ell )  D_\ell$

$$\aligned & \left |   b_{\ell +1 }
  \right |\le  \left |  b_{\ell  }   \right
 |
+  \left |
   \Omega  _{\ell  }  (\omega ,
z,\Phi _{\ell +1} )\right |  \le C_{18}(\ell )  |z|^2  .
 \endaligned \tag 4.18$$
The same bounds hold for $c_{\ell +1}$ and $d_{\ell +1}$.  Consider
now

$$\aligned &\Omega _{\ell +1}( \omega , z ,f _{\ell +1})
=  \Omega _{\ell  }( \omega , z ,f _{\ell +1 })
\\& +\int _{[0,1]^2}\Omega _{\ell  }''( \omega , z ,sf _{\ell +1 }+t\Phi _{\ell +1}( \omega , z ))dt\, ds\, \cdot f _{\ell +1 }\cdot
\Phi _{\ell +1}( \omega , z ).
\endaligned \tag 4.19
$$
By induction, by  Cauchy integral formula, by (4.11), by (4.16) and
by (4.5) we have
$$ \aligned &  | \Omega  _{\ell  }''( \omega , z ,sf _{\ell +1 }+t\Phi _{\ell +1}( \omega , z ))  \cdot f _{\ell +1 }\cdot \Phi
_{\ell +1}( \omega , z )|\le \\& \le C_A(\ell ) C_\Omega (\ell
)D_\ell  \frac{|z|^{2} \|
  f _{\ell +1 }\|  _{\ell ^{2,-2}} }{(
  e^{-2}2^{-\ell -2}\varepsilon _0)^2 } (2|z|+\|
  f _{\ell +1 }\|  _{\ell ^{2,-2}} )^2 .
\endaligned \tag 4.20
$$
Notice we have used $e^{-2}e^{ 2^{-\ell   } }\varepsilon
_0-e^{-2}e^{ 2^{-\ell   }   -2^{-\ell  -2} }\varepsilon
_0>\varepsilon _0e^{-2}2^{-\ell -2}$ and

$$\| sf _{\ell +1 }+t\Phi _{\ell +1} \|  _{\ell
^{2,-2}} \le \varepsilon _0 e^{-2}2^{-\ell -2 }+e^{-(2-2^{-\ell
-1})}\varepsilon _0<  e^{-2}e^{ 2^{-\ell   }   -2^{-\ell  -2}
}\varepsilon _0.$$
For $ C_{21}(\ell):=   C_\Omega (\ell )  +  e^4
3^4C_A(\ell ) C_\Omega (\ell ) 2^{2\ell+2}D_\ell \varepsilon _0$
$$\aligned & |\Omega  _{\ell +1  }(f_{\ell +1  })|\le
|\Omega  _{ \ell   }(f_{\ell +1  } )| +  e^4 3^4 C_A(\ell ) C_\Omega
(\ell )   2^{2\ell+2}D_\ell \varepsilon _0  |z|  \|
  f _{\ell +1 }\|  _{\ell ^{2,-2}}   \\& \le C_{21}(\ell)
  \left ( |z|+  \| f_{\ell +1  }  \|
_{\ell ^{2, -2}}\right ) \| f_{\ell +1  }  \| _{\ell ^{2, -2}} .
\endaligned \tag 4.21
$$
The same bound holds for $ \max \{ |\Gamma _{\ell +1 }|,|Z _{\ell +1
}|, \| \widehat{\Cal F } _{\ell }\| _{\ell ^{2,2}}\}$ by the same
argument   using

$$ \aligned & \|\widehat{\Cal F } _{\ell  }''( \omega , z ,sf _{\ell
}+t\Phi _{\ell +1}( \omega , z ))  \cdot f _{\ell +1 }\cdot \Phi
_{\ell +1}( \omega , z )\| _{\ell ^{2,2}}\le \text{rhs}(4.20).
\endaligned \tag 4.22
$$
 Finally we consider $A_{\ell
+1}(\omega , z)$. We bound each of the terms in the rhs in (4.14).
By Lemma 4.2  for $\Phi _{\ell +1}$  and (4.18) for $c _{\ell +1}$
$$\aligned & \|c_{\ell +1}
 \sigma _3\Phi _{\ell +1}  \| _{\ell ^{2,2}}\le
 C_{23}(\ell )
\varepsilon _0 ^{-\ell +1   } |z|^{\ell +3 }, \,   C_{23}(\ell ):=
C_{18}(\ell )  C_A(\ell )  D_\ell e^{ 2(\ell -1) } .\endaligned \tag
4.23$$ Similarly
$$ \aligned &\|\widehat{\Cal F}_{\ell } (\omega , z, \Phi _{\ell +1}) \| _{\ell
^{2,2}}\le C_{24}(\ell )   \varepsilon _0 ^{-\ell  +1 } |z|^{\ell +2
}, \, C_{24}(\ell ):= 2C_\Omega (\ell )C_A(\ell ) D_\ell
 e^{2 (\ell
-1) },
\endaligned \tag 4.24$$

$$\aligned &  \|\Phi _{\ell +1}'(\omega ,z ) \cdot
 (b_{\ell+1}, c_{\ell+1} )\| _{\ell ^{2,2}}  \le  \|
 \Phi _{\ell +1}'(\omega ,z )\| _{\ell ^{2,2}}
 (|b_{\ell+1}|+| c_{\ell+1}| ) \\&  \le  8C_{23}  (\ell )
 \varepsilon _0 ^{-\ell  +1 }|z|^{\ell +2
}   .
\endaligned \tag 4.25$$
We have for $|z|\le e^{-2+2^{-\ell -1}} \varepsilon _0$

$$ \aligned &  \|  \widetilde{A} _{\ell
 }(\omega ,z)  \| _{\ell ^{2,2}}\le    \left \|\sum _{m+n\ge \ell +2 }
\frac{\partial _z^m\partial _{\overline{z}}^nA _{\ell }(\omega ,0)
}{m!n!} z^m \overline{z}^n\right \| _{\ell ^{2,2}} \le \\& \le
  \sum _{j\ge \ell +2 } \frac{2^j  |z|^{j}}{j! \left ( \varepsilon
_0 \exp (2^{-\ell}-2)  \right )^{j }}  C_A (\ell ) {\left (
\varepsilon _0 \exp (2^{-\ell}-2)   \right ) ^{2}}\\& \le C_A(\ell )
  \sum _{j\ge \ell +2 } \frac{2^j  |z|^{j}}{j! \left ( \varepsilon
_0 \exp (2^{-\ell}-2)  \right )^{j -2}}     \le  C_{26}(\ell )
\varepsilon _0   ^{-\ell }  |z|^{\ell +2} ,  \\&   C_{26}(\ell
):=C_A(\ell )e^2
    \frac{2^{\ell +2}     }{(\ell +2)!
      \left ( \exp (2^{-\ell}-2) \right ) ^{\ell }  }, \, C_{26}(1
):=
    \frac{c(1) e^{5/2}2^{3}     }{3!
        }
       .\endaligned \tag 4.26
$$
To close the inequalities we need $C_A $, $C_b $, $C_\Omega $
 in $\ell ^\infty (\Bbb N)$ such that
$$\aligned &
\varepsilon _0  e^2 \| C_A(j)  2^{ j +2}D_j \| _{\ell ^\infty }<1\,
, \quad  \varepsilon _0\| C_A(j )D_j e^{
 (2^{-j  -1} -2 )} \| _{\ell ^\infty }  < 1,\\&
  C_b(\ell +1 )\ge C_b (\ell )+
\varepsilon _0  C_\Omega (\ell ) C_A (\ell )  D_\ell ,\\&  C_\Omega
(\ell +1 )\ge   C_\Omega (\ell )  + \varepsilon _0 C_A(\ell )
C_\Omega (\ell ) 2^{2\ell+2}e^4 3^2D_\ell , \\&
  C_A(\ell +1 )\ge C_A(\ell )
    \frac{2^{\ell +2}  e^{2\ell +2}   }{(\ell +2)!
         }+
\\& +8\varepsilon _0^2\left [ C_b (\ell )+
\varepsilon _0  C_\Omega (\ell ) C_A (\ell )  D_\ell \right
]C_A(\ell ) D_\ell e^{ 2(\ell -1) }+\\& +\varepsilon _0C_\Omega
(\ell )C_A(\ell ) D_\ell
 e^{2 (\ell
-1) }
 .
  \endaligned \tag 4.27
$$
We set $C_A(1)=   C_b(1)=C_\Omega(1)=  c(1)$  and   define three
sequences recursively using equalities in the last three
inequalities in (4.27). They are in $\ell ^\infty (\Bbb N)$ because
of the following elementary fact:

\proclaim{ Lemma 4.5} Consider a sequence  $x_{n+1}=e^{x_nd_n
\varepsilon _0}(x_n+ x_n d_n+\varepsilon _0x_n ^N d_n)$ for fixed
$N>1$ and $d \in \ell ^1 (\Bbb N)$. Suppose we have inequalities
$2\varepsilon _0 x_1\| d \| _{1 }\exp ( 2\| d \| _{1} ) <\log 2$ and
$\varepsilon _0 2^{N-1} e^{2(N-1) \| d \| _{1}  } x_1^{N-1}
 <1$. Then $x_n\le 2x_1\exp (2\| d \| _{1} ) $ for
         all $n $. \endproclaim {\it Proof.}
By induction
         $$\aligned &x_{n+1}\le  e^{2\varepsilon
_0 x_1d_n \exp (2 \| d \| _{1} )   } \left ( 1+d_n +\varepsilon
_0d_n 2^{N-1} e^{2(N-1) \| d \| _{1}  } x_1^{N-1} \right ) x_n
\\&  \le e^{2\varepsilon
_0 x_1\| d \| _{1 }\exp (2\| d \| _{1} )+2\| d \| _{1}}x_1\le
2x_1\exp (2\| d \| _{1} ).
 \endaligned $$
               This yields Lemma 4.5.
                Since  we can bound from above our three
                sequences by a sequence satisfying Lemma 4.5, this also
 concludes Lemma 4.3.

 \bigskip \noindent We define $g(t)=f(t)-\sum _{\ell \ge 2}\Phi _{\ell}
 (\omega (t),z(t))$. By
  $|\omega (t)-\omega _0|+|z(t)|\le C \epsilon \ll  e^{-2}\varepsilon
   _0$ the series converges. Furthermore $\partial _tg(t)=\partial _t
   f(t)-\sum _{\ell \ge 2}\partial _t\left ( \Phi _{\ell}
 (\omega (t),z(t))\right ).$
 We have $g\in \ell ^2_c(\Cal H_\omega )$ with
 $\| g \| _{\ell ^2}\lesssim \| f \| _{\ell ^2}+|z|^2$ and
$\| g \| _{\ell ^{2,-2}}\lesssim \| f \| _{\ell ^{2,-2}}+|z|^2$. We
have particular $f_\ell \to g $, $\Cal H_\omega f_\ell \to \Cal
H_\omega  g $  (notice that $\Cal H_\omega $ is a bounded operator),
 and  $\partial _tf_\ell \to \partial _tg$ in $\ell ^2$.
 We have $A_\ell (\omega ,z)\to 0$ uniformly for
$|\omega  -\omega _0|+|z |\le    C\epsilon  $. We set $\Cal F _\ell
( f_\ell )=$
$$\aligned &   =\Cal F   ( g ) + \widehat{\Cal F} _\ell ( g )+ \left (
\Cal F  ( f_\ell ) -\Cal F   ( g )\right ) + \left ( \widehat{\Cal
F} _\ell  ( f_\ell ) -\widehat{\Cal F} _\ell   ( g )\right ) \to
\Cal F   ( g ) + \lim _{\ell \to \infty }\widehat{\Cal F} _\ell ( g
).\endaligned$$   We have

$$\aligned & \widehat{\Cal F} _\ell (g   ) =
\widehat{\Cal F} _1 (g   )+\sum _{j=2}^{\ell }\left ( \widehat{\Cal
F} _j (g   ) -\widehat{\Cal F} _{j-1} (g   ) \right ) .\endaligned$$
Since for a fixed $C$ by (4.19) with $\Omega $ replaced by
$\widehat{\Cal F}$ and by (4.20)
$$\| \widehat{\Cal
F} _j (g   ) -\widehat{\Cal F} _{j-1} (g   ) \| _{\ell ^{2,2}}\le C
D_j 2^{2j +2}e^{4}      |z|^2 \| g \| _{\ell ^{2,-2}}\to 0 \text{
for $j\nearrow \infty$,}$$ we have uniform convergence
$\widehat{\Cal F} _\ell ( g )\to \widehat{\Cal F}  ( g )$ in $\ell
^{2,2}$. If $| z| _{L^\infty _t}\le \varepsilon _0 $ (recall that $|
z| _{L^\infty _t}\le C \epsilon$) we have
 $\| \widehat{\Cal F} (\omega , z,  g )\|
_{\ell ^{2,2}}\le C  |z| \| g\| _{\ell ^{2,-2}} $.  By orbital
stability we know $|z|+\| f \| _{\ell ^2}\lesssim \epsilon$. Since
$\| f-g\| _{\ell ^2} \lesssim |z|^2$ we get $\| g \| _{\ell
^2}\lesssim \epsilon$. Taking the limit for $\ell \nearrow \infty $
in (4.3) we obtain (4.1). Notice that since $\sigma _1f_\ell
=\overline{f_\ell}$ for all $\ell$, we have $\sigma _1g
=\overline{g}$. This ends the proof of Lemma 4.1.

 \head \S 5  $g$ is asymptotically free\endhead

Our first result is the following one:

\proclaim{Lemma 5.1} Consider the function $g$ in Lemma 4.1. We have
for a fixed $C_0$
 $$\aligned &
 \| g  \|    _{\ell ^{\frac{3}{2}r }(
\Bbb {Z}  ,L^{\infty }_{t}([n,n+1],\ell ^{ p} ))}\le C_0 \epsilon \,
\text{ for all admissible $(r,p)$} \\& \|
  g  \| _{L^2_t \ell ^{2,-2}} \le C_0 \epsilon .\endaligned
 $$

\endproclaim

{\it Proof.} Set $P_c( \omega )  =P_c(\Cal H_\omega )  $. Let
$\omega _0$ be as in Theorem  1.2. Set $\vartheta =\omega -\omega
_0+\dot \gamma$. We have
$$\aligned & i  \dot g - (\Cal H_ {\omega _0}+\vartheta  P_c(
 \omega _0)  \sigma _3) g = \vartheta \left [ P_c(
 \omega  )-P_c(
 \omega _0)  \right ]\sigma _3 g + \widehat{\Cal F} (\omega , z,g)+O(|g|^7)
 \endaligned  $$
    We split $P_c( \omega _0 )=P_+(\omega _0)+P_-(\omega _0 )$
with $P_\pm (\omega _0)$ the projections in $\sigma _e(\Cal
H_{\omega _0}) \cap \Bbb R_\pm $, see \S 11. Then we rewrite the
above equation as
$$\aligned & i  \dot g - (\Cal H_ {\omega _0}+\vartheta (P_+(
 \omega _0)-P_-(\omega _0)   ) g =    \vartheta \left [ P_c(
 \omega _0)  \sigma _3 -P_+(
 \omega _0)+P_-(\omega _0)\right ] g\\& +
 \vartheta \left [ P_c(
 \omega  )-P_c(
 \omega _0)  \right ]\sigma _3 g+\widehat{\Cal F} (\omega , z,g)+O(|g|^7).
 \endaligned  $$
We prove in Lemma 11.1 that for any pair $s_1,s_2\in \Bbb R$ there
is $ c_{s_1,s_2} (\omega )$ upper semicontinuous in $\omega$  such
that for  $j=0,1$

$$   \|   P_c (\omega
) \sigma _3-  (P_+ (\omega )-P_- (\omega )) \| _{B (\ell ^{2,s_1}
,\ell  ^{2,s_2}  ) } \le c_{s_1,s_2} (\omega )<\infty . $$ By
orbital stability we have $\vartheta =O(\epsilon )$. Then we can
write the above equation  as
$$\aligned & i  \dot g - (\Cal H_ {\omega _0}+\vartheta (P_+(
 \omega _0)-P_-(\omega _0)   ) g =   \epsilon O_{loc}(g)+O(|g|^7)
 \endaligned  $$
for   $  O_{loc}(g)$ such that $\| O_{loc}(g)\| _{\ell ^{2,2  } }
\le C   \|   g \| _{\ell ^{2,-2 } }$, for some $C$. We set $$ U_\pm
(t,t')= e^{-i(t-t')\Cal H _{\omega_0}} e^{\pm i\int _{t'}^t d\tau
(\dot \gamma (\tau )+\omega (\tau )-\omega _0) }P _{\pm}(\omega _0)
\tag 5.1$$ and we write
$$\aligned & P _{\pm}(\omega _0)g(t)=  U_\pm (t,0)g(0)+\int _0^tU_\pm (t,t')
(\epsilon O_{loc}(g)+O (|g|^7)) dt'.
\endaligned $$
By a standard continuation argument Lemma 5.1 is a consequence of
the following result:

\proclaim{Lemma 5.2} There is a fixed $C$ such that if for all
admissible $(r,p)$
 $$\aligned &
 \| g  \|    _{\ell ^{\frac{3}{2}r }(
\Bbb {Z}  ,L^{\infty }_{t}([n,n+1],\ell ^{ p} ))}\le 2C \epsilon \,
 \, , \quad  \|
  g  \| _{L^2_t \ell ^{1,-2}  } \le 2C \epsilon  \endaligned
  \tag 5.2
 $$
then
$$\aligned &
 \| g  \|    _{\ell ^{\frac{3}{2}r }(
\Bbb {Z}  ,L^{\infty }_{t}([n,n+1],\ell ^{ p} ))}\le  C \epsilon \,
, \quad  \|
  g  \| _{L^2_t \ell ^{1,-2}} \le  C \epsilon  \endaligned
  \tag 5.3
 $$
\endproclaim
{\it Proof.} Let $\| g\| _{(r,p)}=\| g  \|    _{\ell ^{\frac{3}{2}r
}( \Bbb {Z}  ,L^{\infty }_{t}([n,n+1],\ell ^{ p} ))}$. For fixed
$C_{p,\sigma }$
$$\| P_ c (\omega  ) -P_ c (\omega _0  )\| _{B(\ell ^{p,\sigma},
 \ell ^{p,\sigma}
)}\le C_{p,\sigma } |\omega - \omega _0|  . $$ So from $C_{p,\sigma
} |\omega (t) - \omega _0|  \le 1/2$ for all $t$ and  $g (t)=P_ c
(\omega (t) )g(t)$ we conclude $ \| g( t)\| _{\ell ^{p,\sigma }}/ \|
P_c(\omega _0) g( t)\| _{\ell ^{p,\sigma }} \in [1/2,2]$.  From this
we conclude  $ \| g \| _{ {(r,p)}}/  \| P_c(\omega _0) g \| _{
{(r,p)}}$ and $\| g \| _{L^2_t\ell ^{2,-2}}/
 \| P_c(\omega _0) g \| _{L^2_t\ell ^{2,-2}} $ are in $[1/2,2]$.
Hence it is enough to prove $(5.2)\Rightarrow (5.3)$ with $g$
replaced by $ P_c(\omega _0) g$. By Lemmas 3.1-2 we have for fixed
constants  $$\| U_\pm (t,0)g(0)\| _{  (r,p)} + \| U_\pm (t,0)g(0)\|
_{ L^2_t\ell ^{2,-2}} \lesssim \| g (0)\| _{\ell ^2}\lesssim
\epsilon .$$ We have by Lemmas 3.5 and 3.6
$$\aligned & \left | \int _0^tU_\pm (t,t')
\epsilon O_{loc}(g)  dt'\right \| _{(r,p)\cap L^2_t\ell ^{2,-2}}\le
C\epsilon \| g \| _{  L^2_t\ell ^{2,-2}}=O(\epsilon ^2).
\endaligned$$ By Lemma 3.1, $$ \aligned & \left  \| \int _0^tU_\pm (t,t')
  O (|g|^7)  dt'\right \| _{(r,p) } \le   C_0
  \| g^7\| _{L^1_t\ell ^2} \le C_0\| g
\| _{ L^{ \infty }_{t} \ell  ^{2  }  }
 \| g   \| ^6 _{ L^{6}_{t}\ell
^{\infty }  } \\&  \le C_0 \| g \| _{ L^{ \infty }_{t} \ell  ^{2  }
} \| g\| ^6 _{ \ell ^6 (\Bbb {Z},L^\infty _t[n, n+1]),\ell ^\infty
))} = O( \epsilon ^7).
\endaligned  $$
By (b) Lemma 3.3
$$\aligned & \left  \| \int _0^tU_\pm (t,t')
  O (|g|^7)  dt'\right \|  _{  L^2_t\ell  ^{2,-2} }
\le   C  \int _0^\infty \| O(g ^7) (s) \| _{ \ell ^2}ds = O(
\epsilon ^7).
\endaligned  $$
This yields Lemma 5.2 and concludes the proof of Lemma 5.1.

\proclaim{Lemma 5.3} We have $^tg (t)=(h(t),\overline{h}(t))$. There
exists  $r_+ \in \ell ^2$ such that $\|  r_+ \| _{\ell ^2} \le C
\epsilon $ for fixed $C=C(\omega _0)$ and
$$\lim _{t\to  \infty }\| h(t)- e^{ it \Delta}r_{+}\| _{\ell ^2}=0.$$
\endproclaim
{\it Proof.} First of all, all transformations in \S 4 preserve the
symmetry $\sigma _1f_\ell =\overline{f}_\ell .$ Hence also $\sigma
_1g =\overline{g}  .$  For $  U_\pm (t,t') $ defined in (5.1) we
have for
 $t_1<t_2$
$$\aligned & \| U(0,t_2)  g(t_2)-
U(0,t_1)  g (t_1) \| _{\ell ^2}\le  \\& \le \| U_\pm (0,t')
(\epsilon O_{loc}(g)+O (|g|^7))\| _{L^1((t_1,t_2),\ell ^2)} \lesssim
\\&   \| g\| _{ L^2((t_1,t_2), \ell ^{2,-2})} +\| g\| _{
L^6((t_1,t_2), \ell  ^{\infty } )}  \to 0
  \text{ for $t_1\to \infty$}.
\endaligned $$
Then   consider $w_+=P _{+}(\omega _0)w+P _{-}(\omega _0)w$  with
 $P _{\pm}(\omega _0)w=\lim _{t\to \infty}U_\pm (0,t )
   g(t ) .$ We have

$$ \lim _{t \to \infty } \| P_c(\omega _0)g (t) \|_{ \ell ^{2,-2}}
=\lim _{t \to \infty } \| ( U_+(t,0)+U_-(t,0))P_c(\omega _0)g (0)
\|_{ \ell ^{2,-2}}=0
$$
with the second equality true for any $g(0) \in \ell ^2$.  By $ \|
P_c(\omega _0)g (t) \|_{ \ell ^{2,-2}}/ \|  g (t) \|_{ \ell ^{2,-2}}
\in [1/2,2]$ it follows
$$ \lim _{t \to \infty } \|  g (t) \|_{ \ell ^{2,-2}}
= 0
$$
 Combining the above we have for $\theta (t)=\int
_{0}^t d\tau (\dot \gamma (\tau )+\omega (\tau ) )$

$$\lim_{t\to  \infty}\|g (t)-e^{i\left [ (t\omega
_0-\theta (t)+\theta (0) \right ] (P_+(\omega _0)- P_-(\omega _0))}
e^{-it\Cal H _{\omega _0}} w _+\|_{\ell ^2}=0 .$$ We claim that  the
following strong limit exists $$ W(\omega_0)=\lim_{t \nearrow
\infty}e^{ it\Cal H _{\omega_0}}e^{ -it( -\Delta  +\omega
_0)\sigma_3}.\tag 5.4$$ The existence of the above limit follows
from the existence of the strong limits
$$ W_1=\lim_{t \nearrow
\infty}e^{ it\Cal H _{\omega_0}}e^{ -it\Cal H _0} P_c(\Cal H _0)
\text{ and } W_2=\lim_{t \nearrow \infty}e^{  it  \Cal H _0 }e^{ it(
 \Delta  -\omega _0)\sigma_3} \tag 5.5$$ with $\Cal H _0= ( H
+\omega _0)\sigma_3$. The first limit in (5.5) exists by Lemma 9.2.
The second  limit in (5.5)  exists by Pearson's Theorem, see Theorem
XI.7\cite{RS}, from the fact that $H+\Delta =q$ is trace class. Then
the limit in (5.4) exists with $W(\omega_0)=W_1\circ W_2$ by the
"chain rule", Proposition 2 ch.XI \cite{RS}. Furthermore
$W(\omega_0)$ is an isomorphism from $\ell ^2$ to $\ell ^2_c(\Cal H
_0)$. Set $ R_+=W(\omega_0)^{-1}e^{i\theta(0)(P_+(\omega _0)-
P_-(\omega _0))}w_+  .$  Notice that since $e^{ it\omega
_0\sigma_3}$ is a unitary matrix periodic in $t$ and $e^{i t\omega
_0\sigma_3}R_+ $ describes a circle in $\ell ^2$. Then we have
$$\aligned & \lim_{t\to +\infty}
 \|e^{-it\Cal H _{\omega _0}} W(\omega_0)e^{i t\omega
_0\sigma_3}R_+ -e^{ -it( -\Delta  +\omega _0)\sigma_3}e^{i t\omega
_0\sigma_3}R_+\|_{\ell ^2}=0.
\endaligned $$
Since $W(\omega _0)$ conjugates $\Cal H_{\omega _0}$ into $ \sigma
_3(-\Delta +\omega _0)$,  we get
$$e^{( it\omega
_0+i\theta (0))(P_+(\omega _0)- P_-(\omega _0))} e^{-it\Cal H
_{\omega _0}} w_+ =e^{-it\Cal H _{\omega _0}} W(\omega_0)e^{
it\omega _0\sigma_3}R_+ .$$ Hence the last two limits  and the
definition of $R_+$ imply the limit
$$  \lim_{t\to +\infty} \left \|  e^{i\theta (t) \sigma _3}g (t) -
 e^{  it \Delta   \sigma_3}R_+  \right \|_{\ell ^2}=0.$$
Hence expressing the last formula in components we obtain Lemma 5.3.

\bigskip
 \head \S 6  Proof of $|z(t)|\approx |z(0)|$  \endhead
We consider $f=\Phi (\omega , z)+g$. Notice that  $\Phi (\omega ,
z)$ is analytic for $\max \{ |\omega -\omega _0| , |z| \} \le
e^{-2}  \varepsilon _0 .$ The equation for $z$ can be written as
$$ i\dot z -\lambda (\omega )z=a(\omega ,z )+ Z(\omega ,z
,g)   $$  with $|Z(\omega ,z ,g)|\le C |z|\| g\| _{\ell ^{2,-2}}.$
We
write this equation in the form

$$i\dot z -\lambda (\omega )z=a(\omega ,z )+\langle g ,
    A(\omega ,z
 )\rangle  + \langle g ,
\Cal A (\omega ,z,g
 )g\rangle   \tag 6.1$$
with $\| \Cal A (\omega ,z,g
 ) \| _{B(\ell ^{2,-2},\ell ^{2, 2} )}\le C$.
Our aim is to show: \proclaim{Lemma 6.1} There is a fixed constant
$C$ such that $|z(0)|\ge C \sqrt{\epsilon} \| f (0)\| _{\ell ^2}$
implies  $|z(t)|\approx |z(0)|$ for all $t$.
\endproclaim
Lemma 6.1 is an immediate consequence of the following lemma:

\proclaim{Lemma 6.2} There is a change of variables $z=\zeta +\alpha
(\omega ,\zeta ) +\langle g, B (\omega ,\zeta ) \rangle$ with
$|\alpha (\omega ,\zeta )|\le C| \zeta |^2$ and $\langle g, B
(\omega ,\zeta ) \rangle \le C | \zeta | \| g\| _{\ell ^{2,-2}}$ for
a fixed $C$, such that
$$i\dot \zeta  -\lambda (\omega )\zeta =d(\omega ,|\zeta |^2 )\zeta
 + \langle g ,
C(\omega ,\zeta ,g
 )g\rangle   \tag 6.2$$
with $d(\omega ,|\zeta |^2 )$ real valued and $|\langle g , C(\omega
,\zeta ,g
 )g\rangle |\le C \| g\| _{\ell ^{2,-2}}^2$ for a fixed $C$.

\endproclaim
Assuming Lemma 6.2 we have
$$\frac{d}{dt} | \zeta  (t)| ^2= 2\Im \left [ \langle g ,
C(\omega ,\zeta ,g
 )g\rangle  \overline{\zeta }\right ]$$
and so by Lemma 5.1
$$\left |  | \zeta  (t)| ^2 - | \zeta  (0)| ^2\right |\le C_1
\epsilon \| g\| _{\ell ^{2,-2}}^2\le C_2\epsilon \| g(0)\| _{\ell
^{2 }}^2 .$$ Then $|\zeta (0)|\ge   \sqrt{ 2C_2 \epsilon} \| g (0)\|
_{\ell ^2}$ implies $|\zeta (0)|/\sqrt{2}\le |\zeta (t)|\le \sqrt{3}
|\zeta (0)|$ for all $t$. Since $|z(t)|\approx |\zeta (t)|$, this
concludes Lemma 6.1.

  \head \S 7 Proof of Lemma 6.2  \endhead

Lemma 6.1 is a consequence of Lemmas 7.1 and 7.3 below.
\proclaim{Lemma 7.1} There is a change of variables $z=\varsigma
+\alpha (\omega ,\varsigma  )  $ with $|\alpha (\omega ,\varsigma
)|\le C| \varsigma |^2$   for a fixed $C$, such that
$$i\dot \varsigma  -\lambda (\omega )\varsigma =
d(\omega ,|\varsigma |^2 )\varsigma + \langle g , c(\omega
,\varsigma
 ) \rangle  +
C(\omega ,\varsigma ,g
 )    \tag 7.1$$
with $d(\omega ,|\varsigma |^2 )$ real,     $ c(\omega ,\varsigma )
  \in \ell ^{2,2}$
  and   $|C(\omega ,\varsigma ,g
 ) |\le C \| g\| _{\ell ^{2,-2}}^2$ for a fixed $C$.

\endproclaim

{\it Proof.} We write (6.1) in the form
$$\aligned & i\dot z -\lambda (\omega )z=\sum _{m+n\ge 2}
 a_{mn}(\omega ) z^m \overline{z}^n +  \langle g ,
    A(\omega ,z
 )\rangle  + \langle g ,
\Cal A(\omega ,z,g
 )g\rangle  \\&  i\dot \omega = b(\omega , z)+\langle g ,
    B(\omega ,z
 )\rangle  + \langle g ,
\Cal B(\omega ,z,g
 )g\rangle  .\endaligned \tag 7.2$$
We define inductively, with (7.3) for $\ell =1$ equal to (7.2),
$$\aligned & i\dot z _{\ell} -\lambda (\omega )z_{\ell} =
a_\ell (\omega , z_\ell ) +
    \alpha _{\ell} (\omega ,z_{\ell}
   ,g )  \\&  i\dot \omega = b_{\ell} (\omega , z_{\ell} )+\beta (
     \omega ,z_{\ell}
   ,g )  . \endaligned \tag 7.3$$
Let $\delta _0=  \varepsilon _0^2$, with  $\varepsilon _0>0$ the
constant in \S 4. We assume the following {\bf inductive
hypotheses}.

{\item {(1)}} $a_\ell  $,  $b_\ell  $, $ \alpha _\ell $ and $\beta
_\ell $   are analytic functions  in $(\omega , z ,\overline{z})$
for $\max \{ |\omega -\omega _0|,|z|\}\le e^{-(2-2^{-\ell})} \delta
_0;$ if we consider expansions
$$ a_\ell (\omega , z)=\sum _{ m+n\ge 2}
 a_{\ell mn}(\omega ) z^m  \overline{z}^n \, , \quad b_\ell (\omega , z)=\sum _{ m+n\ge 2}
 b_{\ell mn}(\omega ) z^m  \overline{z}^n$$
then for $\omega$ real,  the coefficients $a_{\ell mn}(\omega )  $
and  $b _{\ell mn}(\omega )  $ are real;

{\item {(2)}}  for $I_\ell $ defined by $(m,n)\in I_\ell $ either if
$m+n\ge \ell +1$ or if $m-n=1$, we have
$$ a_\ell (\omega , z)=\sum _{(m,n)\in I_\ell}
 a_{\ell mn}(\omega ) z^m  \overline{z}^n \, , \quad \widetilde{a}_\ell (\omega , z)
:=\sum _{(m,n)\in I_{\ell +1}}
 a_{\ell mn}(\omega ) z^m   \overline{z}^n ;    \tag 7.4$$

{\item {(3)}} the  following estimates hold :
$$\aligned & |a_\ell (\omega , z_{\ell   })|\le
C_a (\ell ) |z_{\ell   }|^2\, , \quad |b_\ell (\omega , z_{\ell
})|\le C_b (\ell ) |z_{\ell   }|^2
\\&
 \max \{ |\alpha _{\ell} (\omega ,z_{\ell},
 g)  |  ,|\beta  _{\ell} (\omega ,z_{\ell},
 g)  | \} \le
 C_\alpha (\ell ) (|z_\ell | + \|g\| _{\ell ^{2,-2}}) \|g\| _{\ell ^{2,-2}}
  .
\endaligned \tag 7.5$$

 \medskip  We set $z_{\ell  +1  }=z_{\ell
 }+\phi _{\ell +1  }(\omega , z_{\ell    })$ with

$$\phi _{\ell +1  }(\omega ,
z_{\ell    })= \sum _{m+n=\ell +1} \frac{a_{\ell mn}(\omega )
z^m_{\ell     } \overline{z}^n _{\ell    }}{ (m-n-1)\lambda (\omega
)} \text{ sum over $ m-n\neq 1$.}
$$
  Then we get equations (7.3) for $ \ell +1 $
 with, for $\phi '\cdot (a,b)= \partial _\omega \phi a +
 \partial _z \phi b - \partial _{\overline{z}}
  \overline{\phi} \, \overline{b}$
$$\aligned & \text{$  b_{\ell +1} (\omega , z_{\ell +1} )=b_{\ell} (\omega ,
z_{\ell    }  )$, $ \, \beta _{\ell +1} (\omega , z_{\ell +1}
,g)=\beta _{\ell} (\omega , z_{\ell   } ,g ) , $}
\\& a_{\ell +1} (\omega , z _ {\ell +1})= \widetilde{a}_{\ell  }
(\omega , z_ \ell  ) -   \phi _{\ell +1}'(\omega ,z_{\ell    }
)\cdot (b_{\ell+1}(\omega , z_{\ell +1}),a_{\ell } (\omega ,z_{\ell
} )), \\&  \alpha _{\ell +1} (\omega ,z_{\ell +1}
 ,g)  =  \alpha _{\ell  } (\omega ,z_{\ell    }
 ,g)      -  \phi _{\ell +1}'(\omega
,z_{\ell    } )\cdot (\beta _{\ell+1}(\omega , z_{\ell +1},g)
,\alpha _{\ell } (\omega ,z_{\ell    },g  )
    ) .
\endaligned  \tag 7.6$$
The transformation is designed so that the inductive hypothesis (2)
holds for $\ell +1$, by elementary computation. Similarly, (1) for
$\ell +1$ follows by the definition of $\phi _{\ell +1  }(\omega ,
z_{\ell })$, formulas (7.6) and (1) for $\ell$. Now we focus on the
estimates.

\proclaim{Lemma 7.2}
 Estimates  (7.5) hold for a fixed $C$ replacing $C_a(\ell )$, $C_b(\ell )$
 and $C_\alpha (\ell )$ for all $\ell $.

\endproclaim
{\it Proof.} Set $\kappa =\sup \{ |\lambda ^{-1}(\omega )|: |\omega
-\omega _0|\le \varepsilon _0 \}$ and $d_\ell := 2^{
2\ell+2}\kappa/\ell !$. We have
$$\aligned & |\phi _{\ell +1}(\omega ,z_{\ell })|
\le \sum _{m+n=\ell +1}\frac{|\lambda ^{-1}(\omega )| }{m!n!
 }|\partial
_z^m \partial _{\overline{z}}^na_\ell (\omega , 0)|  \, |z_{\ell
}|^{\ell +1}\\& \le \frac{C _7(\ell )|z_\ell |^{\ell +1} }{ \left (
\delta _0 \exp (2^{-\ell}-2 )\right )^{\ell -1} } \le C_{7} |z_\ell
|^{2} \, , \, C_{7}(\ell ):= \frac{C _a(\ell
)2^{\ell+1}\kappa}{(\ell +1) !}.
\endaligned \tag 7.7
$$
By the Cauchy integral formula for $|\omega -\omega _0|\le
e^{2^{-\ell -1}-2}\delta _0 $ and proceeding as for (4.12) and
(4.13)  we have
$$\aligned & |\phi _{\ell +1}'(\omega ,z_{\ell })|\le
      C_8(\ell ) |z_\ell | \, , \, C_8(\ell ):= C_a(\ell ) e^2d_\ell .
\endaligned \tag 7.8
$$
By $ z _{\ell+1}=z _{\ell }+\phi _{\ell +1} (\omega ,z_{\ell }) $ we
get

   $$|z _{\ell+1}-z _{\ell }| \le  {C _a(\ell )d_\ell  }  |z_\ell
|^{2}\le   {C_a(\ell ) d_\ell } \delta _0 |z_\ell |   .\tag 7.9$$ We
assume ${C_a (\ell ) d_\ell } \delta _0 <1/2$ $\forall$ $\ell$, see
(7.16). Then $|z _{\ell+1}-z _{\ell }| \le e^ {2C_a(\ell ) d_\ell
\delta _0} |z _{\ell+1} |$. We have

$$\aligned & | z _{\ell }|
\le |z _{\ell+1} | +|z _{\ell+1}-z _{\ell }| \le e^{-(2-2^{-\ell
-1})} \delta _0+{C _a(\ell )d_\ell  }  e^{-2(2-2^{-\ell -1})} \delta
_0^{2}<e^{-(2-2^{-\ell  })} \delta _0
\endaligned $$
by (4.11) if $C _a(\ell )d_\ell e^{   2^{-\ell  }-2}2^{\ell
+1}\delta _0 <1.$ We will assume this inequality for the moment, see
(7.16). By   $(1-{C_a (j )d_j } \delta _0 )^{-1}=1+2{C_a (j )d_j }
\delta _0 $

$$\aligned &    |b_{\ell +1 } (\omega , z _{\ell +1})|= |b_{\ell    }
(\omega , z _{\ell   } ) |\le C_{10}(\ell )
 |z _{\ell+1}|^2, \, C_{10}(\ell ):=   e^{2
    C_a( \ell )
    d_ \ell  \delta _0 }C_b( \ell )
    ,
  \endaligned \tag 7.10$$

$$\aligned &    |\beta _{\ell +1 } (\omega , z _{\ell +1},g)|=
|\beta _{\ell   } (\omega , z _{\ell  },g ) |\le  C_{11}  (\ell )
(|z _{\ell +1  } | +\| g\| _{\ell ^{2,-2}}) \| g\| _{\ell ^{2,-2}}
\endaligned \tag 7.11$$ for  $   C_{11}(\ell
):=C_\alpha (\ell ) e^{2
    C_a (\ell  )
    d_\ell    \delta _0 }$.
We have
$$ \aligned & \widetilde{a}_\ell (\omega , z_\ell ) =
 {a}_\ell (\omega , z_\ell )-\widehat{a}_\ell (\omega , z_\ell )
 \text{ with } \widehat{a}_\ell (\omega , z  )=\sum _{ m+n=\ell +1, m-n\neq 1}
 a_{\ell mn}(\omega ) z^m    \overline{z}^n  .\endaligned $$
By an analogue of (7.7) we have,by $ 2 C_a(\ell  ) d_\ell \delta _0
<1$,
$$ \aligned &   |\widehat{a}_\ell (\omega , z _\ell )|\le   \frac{C _a
(\ell )  2^{\ell+1} }{(\ell +1) !\left ( \delta _0 \exp (2^{-\ell}-2
)\right )^{\ell -1} } |z_\ell |^{\ell +1}\le \\& \le  \frac{C_a(\ell
) 2^{\ell+1}  }{(\ell +1) !} (1-{C_a(\ell ) d_\ell } \delta _0
)^{-2} |z_{\ell +1} |^{2}\le  \frac{C_a(\ell ) 2^{\ell+3} }{(\ell
+1) !}
 |z_{\ell +1}
|^{2} .\endaligned \tag 7.12$$ By $a_{\ell +1} (\omega , z _ {\ell
+1})= \widetilde{a}_{\ell  } (\omega , z_ \ell  ) - \phi _{\ell
+1}'(\omega ,z_{\ell    } )\cdot (b_{\ell+1}(\omega , z_{\ell
+1}),a_{\ell } (\omega ,z_{\ell } ))$,
$$ \aligned &
|a_{\ell +1} (\omega , z _ {\ell +1})|\le |   {a}_{\ell  } (\omega ,
z_ \ell  ) | +  \frac{C_a(\ell ) 2^{\ell+3}  }{(\ell +1) !}
 |z_{\ell +1}
|^{2}+\\& + | \phi _{\ell +1}'(\omega ,z_{\ell    } )| \left |
(b_{\ell+1}(\omega , z_{\ell +1}),a_{\ell } (\omega ,z_{\ell }
))\right |     \le | {a}_{\ell  } (\omega , z_ \ell  ) | +
C_{13}(\ell )|z_{\ell +1}|^{2}\\& C_{13}(\ell ):= \frac{C_a(\ell )
2^{\ell+3}  }{(\ell +1) !}
   +\widetilde{C}_{13}(\ell ) \text{ for $\ell >1$}\\&
\widetilde{C}_{13}(\ell ):= C_8(\ell ) \delta _0 \left ( C _{10}
(\ell )+4C _{a} (\ell )\right )  .
 \endaligned \tag 7.13$$
Then
$$ \aligned &
|a_{\ell +1} (\omega , z _ {\ell +1})|\le {C}_{14}(\ell  ) |z_{\ell
+1} |^{2}\, ,\,  C_{14}(\ell ):= e^{4C_a (\ell ) d_\ell \delta
_0}C_{a}(\ell )+C_{13}(\ell ) .
\endaligned \tag 7.14$$
  By $\alpha _{\ell +1} (\omega ,z_{\ell +1}
 ,g)  =  \alpha _{\ell  } (\omega ,z_{\ell    }
 ,g)      -  \phi _{\ell +1}'(\omega
,z_{\ell    } )\cdot (\beta _{\ell+1}(\omega , z_{\ell +1},g)
,\alpha _{\ell } (\omega ,z_{\ell    },g  )),$ $$ \aligned &|\alpha
_{\ell +1} (\omega ,z_{\ell +1 } ,g) |\le    C_{15}(\ell ) (|z_{\ell
+1 }|+\|g\| _{\ell ^{2,-2}}  ) \|g\| _{\ell ^{2,-2}} \\& C_{15}(\ell
):=e^{4C_a (\ell ) d_\ell \delta _0}C_{\alpha}(\ell )+ \delta _0
C_a(\ell ) e^2d_\ell (C _{11}(\ell ) +C_{\alpha}(\ell )).\endaligned
\tag 7.15$$

To close the inequalities we need to find sequences $C_a$, $C_b$ and
$C_\alpha$ in $\ell ^\infty (\Bbb N)$ such that for all $\ell$ and
for $d_\ell :=2^{ 2\ell+2}\sup _\omega |\lambda ^{-1}(\omega )|/\ell
!$
$$\aligned & \,C_a(\ell )
d_\ell   \delta _0 <1/2 \, , \, C _a(\ell )d_\ell e^{ 2^{-\ell
 }-2}2^{\ell +1}\delta _0 <1\\& C_a(\ell +1)\ge C_{14}(\ell )\, ,
\\& C_b(\ell +1)\ge C_{10}(\ell ) \, , \,  C_{\alpha} (\ell +1)\ge
\max \{ C_{11}(\ell ), C_{15}(\ell) \} .
\endaligned \tag 7.16$$
For $\delta _0 $ small there are such $C_a$, $C_b$ and $C_\alpha$ in
$\ell ^\infty (\Bbb N)$ satisfying (7.16). For example the sequence
defined   inductively choosing equality in the last two lines in
(7.16) and  by $C_a(1)=C_b(1)=C_\alpha (1)=c(1)$. This follows from
the fact that these sequences are bounded from above by a sequence
satisfying Lemma 4.5. This concludes the proof of Lemma 7.2.

\bigskip
 The $z_\ell $ converge to $\varsigma $. Taking the limit for $\ell \nearrow
 \infty $ in (7.3) and proceeding as at the end of \S 4 we see that
$\varsigma $  satisfies a system of the form
$$i\dot \varsigma  -\lambda (\omega )\varsigma =
d(\omega ,|\varsigma |^2 )\varsigma + \gamma (\omega ,\varsigma
 ,g)    $$
with $|\gamma (\omega ,\varsigma
 ,g)   |\le C |  \varsigma | \,  \| g\|
_{\ell ^{2,2}} $ and analytic   in $(\omega , \varsigma ,
\overline{\varsigma}  , g)$ for $\max \{ |\omega -\omega
_0|,|\varsigma |  \} \le \delta _0^2.$ Expanding $\gamma (\omega
,\varsigma
 ,g)  $ we get (7.1). This concludes the proof of Lemma 7.1.

\bigskip

\proclaim{Lemma 7.3} There is a change of variables $\varsigma
=\zeta +\langle g, B (\omega ,\zeta ) \rangle$ with  $\langle g, B
(\omega ,\zeta ) \rangle \le C | \zeta | \| g\| _{\ell ^{2,-2}}$ for
a fixed $C$, such that
$$i\dot \zeta  -\lambda (\omega )\zeta =d(\omega ,|\zeta |^2 )\zeta
 + C(\omega ,\varsigma ,g
 )   \tag 7.17$$
with $d(\omega ,|\zeta |^2 )$ real valued and $|C(\omega ,\varsigma
,g
 ) |\le C \| g\| _{\ell ^{2,-2}}^2$ for a fixed $C$.

\endproclaim
{\it Proof.} We write (7.3) in the form
$$\aligned & i\dot \varsigma  -\lambda (\omega )\varsigma =
d(\omega ,|\varsigma |^2 )\varsigma + \sum _{m+n\ge 1}\langle g ,
c_{mn}(\omega
 ) \rangle \varsigma ^m \overline{\varsigma} ^n+
\Cal C(\omega ,\varsigma ,g
 )    \\&  i\dot \omega = b(\omega , \varsigma )+
      \Cal B(\omega ,\varsigma ,g
 )    \endaligned \tag 7.18$$
with $| \Cal B(\omega ,\varsigma ,g
 )  |\le C   \| g\| _{\ell ^{2,-2}}^2$ for a fixed $C$.
We define inductively the following system, which for $\ell =1$
coincides with (7.18):
$$\aligned & i\dot \varsigma _{\ell} -
\lambda (\omega )\varsigma _{\ell} = d(\omega ,|\varsigma _{\ell}|^2
)\varsigma _{\ell} +\langle g , C_{\ell
 }(\omega
, \varsigma _\ell ) \rangle +   \Cal C_{\ell}(\omega ,\varsigma
_{\ell}  ,g
 )   \\&  i\dot \omega = b   (\omega , \varsigma
_{\ell})+ \Cal B_{\ell}(\omega ,\varsigma _{\ell},g
 )  \\& i  \dot g -
(\Cal H_\omega +\dot \gamma P_c(\Cal H_\omega )  \sigma _3) g  =
  {  G}  _\ell
(\omega, \varsigma _\ell ) g  + {\Cal G}  _\ell (\omega, \varsigma
_\ell , g )+O(|g|^7) .\endaligned \tag 7.19$$ Let $\delta _1=
 \delta _0^2$, with $\delta _0$ the constant of Lemmas 7.1--2.
We assume the following {\bf inductive hypotheses}.

{\item {(1)}} $  C_\ell (\omega , \varsigma )  $ and  $B_\ell
(\omega , \varsigma ) $ are analytic  in $(\omega , \varsigma ,
\overline{\varsigma } )$ in $\max \{ |\omega -\omega _0|,|\varsigma
|  \} \le e^{-(2-2^{-\ell})} \delta _1 $ with values in $\ell
^{2,2}$;   $ C_\ell (\omega , \varsigma )  $ belongs to $\ell
^2(\Cal H _\omega ^*)$;

{\item {(2)}} the  following estimates
 hold:
$$\aligned & \| C_\ell (\omega , \varsigma _{\ell   })
\| _{\ell ^{2,2}} \le C_C (\ell )e^{ (\ell -1)(2-2^{-\ell})} \delta
_1 ^{-\ell +1 }    |\varsigma _\ell |^\ell  \, , \,
 | \Cal C  _{\ell} (\omega ,\varsigma _{\ell},
 g)  |  \le
 C _{\Cal C}(\ell )  \|g\| _{\ell ^{2,-2}}^2 \, , \\&  |\Cal B  _{\ell} (\omega ,\varsigma _{\ell},
 g)  |  \le
 C _{\Cal B} (\ell ) (|\varsigma _\ell | +\|g\| _{\ell ^{2,-2}})
   \|g\| _{\ell ^{2,-2}} \, ,     \\&  \|   G _\ell
(\omega, \varsigma _\ell ) \| _{B (\ell ^{2, -2}, \ell ^{2, 2}) }
\le C _{  G} (\ell) | \varsigma _\ell |   \, , \,   \|  {\Cal G}
_\ell (\omega, \varsigma _\ell , g ) \| _{\ell ^{2, 2}} \le C _{\Cal
G} (\ell) \|g\| _{\ell ^{2,-2}}^2
  .
\endaligned \tag 7.20$$
These hypotheses hold for $\ell =1$  with $C _{ C}(1)=C _{\Cal
C}(1)=C _{\Cal B}(1 )=C _{  G} (1)  =C _{\Cal G} (1)=c(1) $ for some
constant $c(1)$.
 We expand

$$  C_{\ell
 }(\omega
, \varsigma _{\ell })   =  \sum _{m+n\ge \ell }  c_{\ell mn}(\omega
 )  \varsigma ^m _{\ell } \overline{\varsigma}_{\ell }  ^n \, , \,   \widetilde{C}_{\ell
 }(\omega
, \varsigma _{\ell }):= \sum _{m+n\ge \ell +1 } c_{\ell mn}(\omega
 )   \varsigma ^m _{\ell } \overline{\varsigma}_{\ell }  ^n
. $$ We set $\varsigma _{\ell +1  }=\varsigma _{\ell }+\phi _{\ell
+1 }(\omega , \varsigma _{\ell   },g)$

$$\aligned & \phi _{\ell +1  }(\omega ,
\varsigma _{\ell   },g)= \langle g, \Phi _{\ell +1  }(\omega ,
\varsigma _{\ell   } ) \rangle \text{ with } \Phi _{\ell +1 }(\omega
, \varsigma _{\ell   } )= \sum _{m+n=\ell  }    \varsigma ^m_{\ell }
\overline{\varsigma }^n _{\ell   }  \gamma _{\ell mn}(\omega )
  \\& \gamma _{\ell mn}(\omega ) =R_{\Cal H _\omega
^*}((n-m+1)\lambda (\omega )) c _{\ell mn}(\omega ).
\endaligned
$$
By induction  $C_\ell (\omega , \varsigma _\ell ) \in \ell _c^2(
\Cal H _\omega ^*)$ and so also $ c _{\ell mn}(\omega ) \in \ell
_c^2( \Cal H _\omega ^*)$ for all $(m,n)$. Then   $ \gamma _{\ell
mn}(\omega )$ and $ \Phi _{\ell +1 }(\omega , \varsigma _{\ell   } )
\in \ell ^2_c( \Cal H _\omega ^*)$ .   We get equations (7.19) for
$\ell +1 $
 with:  $$\align     &   \Cal B _{\ell +1}(\omega , \varsigma _{\ell +1},g)
 =\Cal B _{\ell}(\omega , \varsigma _{\ell  } ,g )+   b
 (\omega , \varsigma _ \ell  )-  b
 (\omega , \varsigma _{\ell +1}   ) \, ;\tag 7.21\\&  G  _{\ell +1}
(\omega, \varsigma _{\ell +1}   ) =  G  _\ell (\omega, \varsigma
_\ell ) \, , \quad
  {\Cal G} _{\ell +1} (\omega
  , \varsigma _{\ell +1}  , g)=  {\Cal G} _{\ell  } (\omega
  , \varsigma _\ell  , g)   ; \tag 7.22\\ &
    C_{\ell +1} (\omega ,\varsigma _{\ell +1}
 ) =
    \widetilde{C}_{\ell} (\omega ,\varsigma _{\ell +1 }) -
    [ G  _\ell (\omega, \varsigma
_\ell ) ]^*
  \Phi _{\ell +1} (\omega , \varsigma _\ell )    -\tag 7.23\\&
     \quad   \quad \quad       - \Phi _{\ell +1}' (\omega
,\varsigma _{\ell } ) \cdot (b (\omega , \varsigma _{\ell  } ) ,d
(\omega , |\varsigma _{\ell } |^2 )\varsigma _{\ell   }
  \  )  + d^{(1)}(\omega , \varsigma _{\ell +1  } )  \,  ;
   \\ &
   \Cal  C_{\ell +1} (\omega ,\varsigma _{\ell +1},
 g)  =
    \Cal C_{\ell} (\omega ,\varsigma _ \ell  ,g)   +\langle g ,
    \widetilde{C}_{\ell} (\omega ,\varsigma _{\ell   }) -
    \widetilde{C}_{\ell} (\omega ,\varsigma _{\ell +1 }) \rangle
      \tag 7.24\\& - \left \langle g,  \Phi _{\ell +1}'
(\omega ,\varsigma _{\ell   }    )\cdot  \left (\Cal B _{\ell
 }(\omega , \varsigma _{\ell  },g),  \langle g, C _{\ell }(\omega
, \varsigma _{\ell  })\rangle +   \Cal C_{\ell } (\omega ,\varsigma
_ \ell  ,g ) \right )    \right \rangle    \\& -\langle  {\Cal G}
_\ell (\omega
  , \varsigma _\ell  , g)  +O(|g|^7),  \Phi
_{\ell +1} (\omega , \varsigma _\ell )\rangle
  +D^{(2)}(\omega
,\varsigma _{\ell +1 },g )  \endalign $$ with   $d^{(1)}(\omega ,
\varsigma _{\ell +1  } )=
\partial _g | _{g=0}  \left [ d(\omega ,|\varsigma _{\ell}|^2
)\varsigma _{\ell}-d(\omega ,|\varsigma _{\ell +1}|^2 )\varsigma
_{\ell +1}\right ] $  and $$ \aligned &     D^{(2)}(\omega
,\varsigma _{\ell +1 },g ) =
 d(\omega
,|\varsigma _{\ell}|^2 )\varsigma _{\ell}-d(\omega ,|\varsigma
_{\ell +1}|^2 )\varsigma _{\ell +1}   -\left \langle
g,d^{(1)}(\omega , \varsigma _{\ell +1  } ) \right \rangle .
\endaligned $$
Assumption (1) for $\ell +1$ holds by an argument  analogous to
Lemma 4.3. We have
$$\aligned & |\phi _{\ell +1}(\omega ,\varsigma _\ell ,g)|\le \sum _{m+n= \ell }
\frac{C_{4.9} }{m!n!
 }\|\partial
_{\varsigma _\ell} ^m \partial _{\overline{\varsigma _\ell }}^nC
_{\ell  } (\omega , 0)\| _{\ell ^{2,2}} \, |\varsigma _\ell |^{\ell
}\| g\| _{\ell ^{2,-2}}  \\& \le C_{4.9}|\varsigma _\ell |^{\ell }\|
g\| _{\ell ^{2,-2}}    \frac{2^\ell  }{\ell !
 }   (e^{ 2^{-\ell }-2}\delta _1) ^{-\ell }  \sup _{  |z|
   \le
 e^{ 2^{-\ell }-2}\delta _1}  \| C _{\ell   } (\omega , z)\| _{\ell ^{2,2}}
  \\& \le \frac{C_C(\ell )C_{4.9}
2^{\ell }  }{\ell ! }(e^{ 2^{-\ell }-2}\delta _1) ^{-\ell +1}
|\varsigma _\ell |^{\ell }\| g\| _{\ell ^{2,-2}} \text{ for
$|\varsigma _\ell |\le e^{ 2^{-\ell }-2}\delta _1$}.
\endaligned
$$
By a similar argument for  $|\varsigma _\ell |\le e^{ 2^{-\ell
}-2}\delta _1$
$$\aligned & | \partial _{(\varsigma , \overline{\varsigma} )}
\phi _{\ell +1}(\omega ,\varsigma _\ell ,g)|\le  \frac{C_C(\ell
)C_{4.9} 2^{\ell +2} }{(\ell -1) ! }(e^{ 2^{-\ell }-2}\delta _1)
^{-\ell +1} |\varsigma _\ell |^{\ell -1}\| g\| _{\ell ^{2,-2}}
\endaligned
$$
 and for $|\varsigma _\ell |\le e^{ 2^{-\ell }-2}\delta _1$ and
$|\omega -\omega _0 |\le e^{ 2^{-\ell -1 }-2}\delta _1$, by an
inequality similar to (4.11),
$$\aligned & | \partial _{\omega}
\phi _{\ell +1}(\omega ,\varsigma _\ell ,g)|\le
   \frac{C_C(\ell )C_{4.9} 2^{ \ell }  e^{2^{-\ell -1 }} }{\ell
! e^{-2}\delta _1 } |\varsigma _\ell |^{\ell   }\| g\| _{\ell
^{2,-2}} .
\endaligned
$$
 Then, for $|\varsigma _\ell |\le e^{ 2^{-\ell }-2}\delta
_1$, $|\omega -\omega _0 |\le e^{ 2^{-\ell -1 }-2}\delta _1$
 and $C_{25}(\ell )
:= {C_C(\ell )C_{4.9} 2^{\ell +1}  }/{\ell ! }$

$$\aligned  & |\phi _{\ell +1}(\omega ,\varsigma _\ell ,g)|\le
C_{25}(\ell ) (e^{ 2^{-\ell }-2}\delta _1) ^{-\ell +1} |\varsigma
_\ell |^{\ell }\| g\| _{\ell ^{2,-2}} \le C_{25}(\ell )  |\varsigma
_\ell | \| g\| _{\ell ^{2,-2}}  \\ & |  \phi _{\ell +1} '(\omega
,\varsigma _\ell ,g)|\le 2^{12}\ell C_{25}(\ell )(e^{ 2^{-\ell
}-2}\delta _1) ^{-\ell +1} |\varsigma _\ell |^{\ell -1}\| g\| _{\ell
^{2,-2}}\le 2^{12} \ell C_{25}(\ell )    \| g\| _{\ell ^{2,-2}} .
\endaligned \tag 7.25 $$
We have
 $    \varsigma _{\ell +1}- \varsigma _{\ell
 } =      \phi _{\ell
 +1} (\omega , \varsigma _{\ell
 },g)$.
So for a fixed $C_0$ from Lemma 5.1

$$\aligned &|\varsigma  _{\ell+1}-\varsigma  _{\ell
}| \le C_{25}(\ell )
 \| g\| _{\ell ^{2,-2}} |\varsigma  _{\ell }|\le
C_{26}(\ell )
 \epsilon |\varsigma  _{\ell }| \\& \le
 e ^{2C_C(\ell )d_\ell \epsilon}
 C_C(\ell )d_\ell \epsilon |\varsigma  _{\ell +1 }|\, , \, C_{26}(\ell ):=
 C_0C_{25}(\ell ) \, , \,
 d_\ell :=  {C_0C_{4.9} 2^{\ell } }/{ \ell   ! }, \endaligned \tag 7.26$$
 where we   assume $2 C_{26}(\ell )
  \epsilon =C_C(\ell ) d_\ell \epsilon < 1,$   see (7.34).
We get  for $h(\omega , \varsigma )=b(\omega , \varsigma  ),
d(\omega , |\varsigma |)\varsigma$

$$\aligned & |h(\omega , \varsigma _{\ell+1})-h(\omega , \varsigma
_{\ell })|\le \| D_\varsigma h(\omega ,\varsigma  ) \| _{L^\infty \{
|\varsigma |\lesssim | \varsigma _{\ell+1}| \} } | \varsigma
_{\ell+1}-\varsigma _{\ell } |  . \endaligned  $$ This yields for $
C _{27}(\ell ):=e ^{2C_C(\ell )d_\ell \epsilon} K(b,d)C_0^{-1}
C_C(\ell )d_\ell $ for a fixed   $K(b,d)$

$$\aligned & |b(\omega , \varsigma _{\ell+1})-b(\omega , \varsigma
_{\ell })|\le     C _{27}(\ell ) \delta _1
 \| g\| _{\ell ^{2,-2}}  | \varsigma
_{\ell+1}  |   \, , \,    \\&  \|d^{(1)}(\omega , \varsigma
_{\ell+1}) \|  _{\ell ^{2, 2}}\le  C _{27}(\ell )    \delta _1
    | \varsigma
_{\ell+1}  |   \, , \\& |D^{(2)}(\omega , \varsigma _{\ell+1},g)  |
 \le C _{27}(\ell )    \delta _1 \| g\| _{\ell ^{2,-2}}^2.
\endaligned \tag 7.27$$
  $\Cal B _{\ell +1}(\omega , \varsigma _{\ell +1},g)
 =\Cal B _{\ell}(\omega , \varsigma _{\ell  } ,g )+   b
 (\omega , \varsigma _ \ell  )-  b
 (\omega , \varsigma _{\ell +1}   )$ implies
$$ \aligned &  | \Cal B_{\ell +1 } (\omega , \varsigma _{\ell+1},g)
 |    \le  C _{28}(\ell )(|\varsigma _{\ell +1}
 | +\| g\| _{\ell ^{2,-2}})  \| g\| _{\ell ^{2,-2}}\, ,\\&
  C _{28}(\ell )= e ^{2C_C(\ell )d_\ell
\epsilon} C _{\Cal B}(\ell )  +  C _{27}(\ell ) \delta _1.
\endaligned \tag 7.28
$$
   We bound
$$\aligned &
\| \widetilde{C}_{\ell } (\omega , \varsigma  ) \| _{ \ell ^{2,2}}
\le \sum _{m+n\ge \ell +1} \frac{|\varsigma |^{m+n}}{m! n!} \|
\partial  _{\varsigma }^m \partial  _{\overline{\varsigma} }^n
 {C}_{\ell } (\omega , 0  ) \| _{ \ell ^{2,2}}\le \\ & \le
C_C(\ell )\sum _{j\ge \ell +1} \frac{2^j |\varsigma |^j}{j! \left
(\delta _1\exp (2^{-\ell}-2) \right ) ^{j-1 }}\le \frac{e^2C_C(\ell
)2^{\ell +1}|\varsigma |^ {\ell +1}}{ (\ell +1) !\left (\delta
_1\exp (2^{-\ell}-2) \right ) ^{\ell  }}
 .
\endaligned \tag 7.29
$$
By a similar argument for $C_{30}(\ell )= {e^2C_C(\ell )2^{\ell +2}
}/{ \ell ! }$

$$\aligned &
\|  \partial _{\varsigma , \overline{\varsigma}}\widetilde{C}_{\ell
} (\omega , \varsigma ) \| _{ \ell ^{2,2}} \le  C_{30}(\ell ) \left
(\delta _1\exp (2^{-\ell}-2) \right ) ^{-\ell  } |\varsigma |^ {\ell
}
 .
\endaligned \tag 7.30
$$
By (7.23) for a constant $K=K(b,d)$  we have
$$\aligned  &
  \| C_{\ell+1} (\omega , \varsigma _{\ell+1}) \| _{
\ell ^{2,2}}\le \| \widetilde{C}_{\ell } (\omega , \varsigma
_{\ell+1}) \| _{ \ell ^{2,2}} +
 \|  \Phi _{\ell +1}' (\omega , \varsigma _\ell )\| _{\ell ^{2,2}}
 K  | \varsigma
 _\ell |^2\\& +   \|  [  {  G}  _\ell  (\omega
  , \varsigma _\ell   ) ]^*  \| _{ B(\ell ^{2,2} ,  \ell ^{2,2} )}
 \|  \Phi _{\ell +1} (\omega , \varsigma _\ell )\| _{\ell ^{2,2}} + \|d^{(1)}(\omega , \varsigma
_{\ell+1}) \|  _{\ell ^{2, 2}}
\\&  \le C  _C(\ell )  \left (
\frac{e^2 2^{\ell +1}}{  (\ell +1)!}    |\varsigma _{\ell +1} |   +
2^{14}\ell C_{25}(\ell )K \delta _1  |\varsigma _{\ell +1} | \right
)
\\& + C _{  G}(\ell ) C_{25}(\ell )2^{\ell +1}
 \delta _1  |\varsigma
_{\ell +1} |  + C _{27}(\ell )    \delta _1
    | \varsigma
_{\ell+1}  |  =  C_{31}(\ell)     |\varsigma _{\ell +1} |
\endaligned \tag 7.31$$
with $C_{31}(\ell)$ defined so that the constants match and where we
used (7.26) and the assumption $ C_0C_C(\ell )d_\ell \epsilon <1/2.$
  By
(7.24) we have,
$$\aligned & |\Cal  C_{\ell +1} (\omega ,\varsigma _{\ell +1},
 g) | \le
    |\Cal C_{\ell} (\omega ,\varsigma _ \ell  ,g) |  +
    \| \widetilde{C}_{\ell} (\omega ,\varsigma _{\ell   }) -
    \widetilde{C}_{\ell} (\omega ,\varsigma _{\ell +1 })
    \| _{\ell ^{2,2}} \| g
    \| _{\ell ^{2,-2}}
     + \\& \|  \Phi _{\ell +1}'
(\omega ,\varsigma _{\ell   }    ) \| _{\ell ^{2,2}} \left ( |\Cal B
_{\ell
 }(\omega , \varsigma _{\ell  },g)| +  \|  C _{\ell }(\omega
, \varsigma _{\ell  })\| _{\ell ^{2,2}} \| g
    \| _{\ell ^{2,-2}}    +  | \Cal C_{\ell } (\omega
,\varsigma _ \ell  ,g )| \right )  \| g
    \| _{\ell ^{2,-2}}      \\& +\|  \Phi _{\ell +1}
(\omega ,\varsigma _{\ell   }    ) \| _{\ell ^{2,2}} \| {\Cal G}
_\ell (\omega
  , \varsigma _\ell  , g)  +O(|g|^7) \| _{\ell ^{2,-2}} +
  |D^{(2)}(\omega , \varsigma _{\ell+1},g)  |
   . \endaligned $$
We have
$$\aligned &
\| \widetilde{C}_{\ell} (\omega ,\varsigma _{\ell   }) -
    \widetilde{C}_{\ell} (\omega ,\varsigma _{\ell +1 })
    \| _{\ell ^{2,2}} \le C _{31}(\ell )
    |\varsigma _{\ell +1 }-\varsigma _{\ell   }| \le C _{31}(\ell )
    C _{25}(\ell ) \delta _1\| g\| _{\ell ^{2,-2}}
    ; \\& \|  \Phi _{\ell +1}'
(\omega ,\varsigma _{\ell   }    ) \| _{\ell ^{2,2}} \left ( |\Cal B
_{\ell
 }(\omega , \varsigma _{\ell  },g)| +  \|  C _{\ell }(\omega
, \varsigma _{\ell  })\| _{\ell ^{2,2}} \| g
    \| _{\ell ^{2,-2}}    +  | \Cal C_{\ell } (\omega
,\varsigma _ \ell  ,g )| \right )   \\& \le 2^{15 }  \ell   C
_{25}(\ell ) \delta _1  \left ( C _{\Cal B} (\ell )+C _{C} (\ell )
+C _{\Cal C} (\ell )\right ) \| g
    \| _{\ell ^{2,-2}}  ;  \\&  \| \Phi _{\ell +1}
(\omega ,\varsigma _{\ell   }    ) \| _{\ell ^{2,2}} \| {\Cal G}
_\ell (\omega
  , \varsigma _\ell  , g)  +O(|g|^7) \| _{\ell ^{2,-2}} \le  C
_{25}(\ell ) \delta _1\left ( C _{\Cal G} (\ell ) +c_0 \epsilon
^5\right ) \| g
    \| _{\ell ^{2,-2}} ^2;\\& |D^{(2)}(\omega , \varsigma _{\ell+1},g)  |
 \le C _{27}(\ell )    \delta _1 \| g\| _{\ell ^{2,-2}}^2
   . \endaligned $$
So
$$\aligned &
   |\Cal  C_{\ell +1} (\omega ,\varsigma _{\ell +1},
 g) | \le C  _{32}(\ell )  \| g \|  _{\ell ^{2,-2}}^2  , \,
C  _{32}(\ell ) := C _{\Cal C}(\ell )+C _{31}(\ell )
    C _{25}(\ell ) \delta _1+\\&  2^{15 }  \ell   C
_{25}(\ell ) \delta _1  \left ( C _{\Cal B} (\ell )+C _{C} (\ell )
+C _{\Cal C} (\ell )\right )+ C _{25}(\ell ) \delta _1\left ( C
_{\Cal G} (\ell ) +c_0 \epsilon ^5\right ) +C _{27}(\ell )    \delta
_1.
\endaligned\tag 7.32$$
We have
$$ \aligned & \| G _{\ell +1} (\omega
  , \varsigma _{\ell +1}   ) \| _{B (\ell ^{2, -2}, \ell ^{2,  2}) }
  = \|G _{1  } (\omega
  , \varsigma _1    )\| _{B (\ell ^{2, -2}, \ell ^{2,  2}) }\le  c(1)  | \varsigma _1 |
   \\& \le  C_{33}(\ell )  | \varsigma _{\ell +1} |
   \, , \,  C_{33 }(\ell ):= c(1)e^{2  \epsilon
\| C_C \| _{\infty (j\le \ell )}
 \| d_j \| _{1}}\\&
 \|  {\Cal G} _{\ell +1} (\omega
  , \varsigma _{\ell +1}  , g) \| _{\ell ^{2, 2}} = \|
  {\Cal G} _{1 } (\omega
  , \varsigma _1  , g)\| _{\ell ^{2, 2}}\le  c(1)
 \|g\| _{\ell ^{2,-2}}^2  .
\endaligned \tag 7.33$$

Now we need for all $\ell$
$$\aligned &
  C_C(\ell )  \epsilon d_\ell <1/2 \, , \,   C_C(\ell )d_\ell
 e^2 2^{\ell +1}\epsilon < \delta _1 \\&
 C_{  C}(\ell +1)\ge C_{31}(\ell )   \, , \,
  C_{\Cal B}(\ell +1)\ge C_{28}(\ell )\, , \\&
C_{\Cal C}(\ell +1)\ge C_{33}(\ell ) \, , \, C _{  G} (\ell +1) \ge
C_{33} (\ell )\, , \, C _{\Cal G} (\ell +1) = c(1).
\endaligned \tag 7.34$$
There exist $C_C,$ $C_{\Cal C} $,  $C_{\Cal B}$  and $C _{\Cal G}$
in $\ell ^\infty (\Bbb N)$  satisfying (7.34), see below (7.16).
Then we can replace the constants in (7.20) with a fixed constant
$C$. Then $\varsigma _\ell \to \zeta $ which satisfies the statement
of Lemma 7.3.

\bigskip
 \head \S 8  Proof of Theorem 1.2  \endhead

We still need to discuss the equation  for  $\omega (t)$. Recall
that we have

$$i\dot \omega  =a(\omega ,\zeta )+\langle g ,
    A(\omega ,\zeta
 )\rangle  + \Cal
A(\omega ,\zeta ,g
 )    \tag 8.1$$ with
 $|\Cal A(\omega ,\zeta ,g
 )|\le C \| g \| _{\ell ^{2,2}}  $ and $\|  A(\omega ,\zeta
 ) \| _{\ell ^{2,2}}\le C |\zeta |$, with $\Cal A(\omega ,\zeta ,g
 )$ and $A(\omega ,\zeta
 )$ analytic in $(\omega ,\zeta , \overline{\zeta} ,g
 )$.
The first step is:

\proclaim{Lemma 8.1} There is a change of variables $\omega =\varpi
+\alpha (\varpi  ,\zeta ) +\langle g, B (\varpi ,\zeta ) \rangle$
with $|\alpha (\varpi ,\zeta )|\le C| \zeta |^2$ and $\|  B (\varpi
,\zeta ) \| _{\ell ^{2, 2}} \le C | \zeta | $ for a fixed $C$, such
that
$$i\dot \varpi    =
D(\varpi ,\zeta ,g
 )    \tag 8.2$$
with   $|  D(\varpi ,\zeta ,g
 )  |\le C \| g\| _{\ell ^{2,-2}}^2$ for a fixed $C$.

\endproclaim
An immediate consequence of Lemmas 8.1 and 5.1 is:

 \proclaim{Corollary  8.2} There is a
fixed $C$ such that $ \|
 \dot \varpi   \| _{L^1\cap L^\infty }<C \epsilon ^2.
$

\endproclaim

 We have: \proclaim{Corollary  8.3} For any $\sigma >0$ there is a   fixed $C$ such that for
$|z(0)|\ge \epsilon $ and $\| f (0) \| _{\ell ^2}\lesssim \epsilon $
we have
$$\inf _{\kappa , \mu }\| u(t) -e^{i\kappa} \phi   _{\mu }
\| _{\ell  ^{2,-\sigma }} \ge C \epsilon .$$
\endproclaim
{\it Proof.} We have for $^t\xi =(\xi _1, \xi _2)$
$$u(t)= e^{i\theta (t)}\left (\phi _{\omega (t)} +z\xi _1(\omega (t))+
\overline{z}\xi _2(\omega (t)\right )+ A(\omega (t), z(t))+ h(t) $$
with $\|A(\omega (t), z(t))\| _{\ell ^{2,2}}\le C |z(t)|^2 $ and $
\lim _{t\to \infty }\| h(t)\| _{\ell  ^{2,-\sigma}}=0.$ By $\| \xi
_2\| _{\ell ^{2,2}} \le C  |E_0-\omega _0| $ and $\| \xi _1-\varphi
_1\| _{\ell ^{2,2}} \le C  |E_0-\omega _0| $, Lemma 2.5,
  for
$t\gg 1$ we have $$\| e^{-i\theta (t)}u(t)-  \phi _{\omega (t)}
-z\xi _1(\omega (t))\| _{\ell  ^{2,-\sigma }}\le C\epsilon
|E_0-\omega _0|
 $$ also by $|z| \lesssim \epsilon $, which follows by orbital
 stability, Lemma 2.3. We have
$$ \aligned &\| \phi _{\omega
 } +z\xi _1(\omega  )-e^{i\kappa} \phi   _{\mu }\| _{\ell
 ^{2,-\sigma }}^2
\approx  \| \phi _{\omega
 } +z\xi _1(\omega  )-e^{i\kappa} \phi   _{\mu }\| _{\ell  ^{2 }}^2 \ge \\&
\| \phi _{\omega
 } -  \phi   _{\mu }\| _{\ell  ^{2 }}^2+ |z|^2\|  \xi _1(\omega  )
 \| _{\ell  ^{2 }}^2-2|z| \left |
  \langle    \phi _{\omega
 }-\phi   _{\mu } , \xi _1(\omega  )\rangle \right |.
 \endaligned $$
By $\| \xi _1-\varphi _1\| _{\ell ^{2 }} \lesssim |E_0-\omega _0| $,
by $\langle \varphi _0,\varphi _1\rangle =0$ and by Lemma 1.1, we
have
$$\aligned &
  \langle  \phi _{\omega
 }-\phi   _{\mu }, \xi _1(\omega  )\rangle
 = \left ( (\omega -E_0) ^{\frac{1}{6}} -
 (\mu -E_0) ^{\frac{1}{6}}\right )
 \langle \varphi _0,
 O(\omega -E_0)\rangle \\& + \langle
 O((\omega -E_0)^{\frac{7}{6}})-O((\mu -E_0)^{\frac{7}{6}})
  , \xi _1(\omega  )\rangle = O \left  ( (\omega -\mu )
  (\omega -E_0)^{\frac{1}{6}}\right ) .
\endaligned $$
Hence $\| \phi _{\omega
 } +z\xi _1(\omega  )-e^{i\kappa} \phi   _{\mu }\| _{\ell  ^{2
 }}^2 \gtrsim | \omega -\mu |^2 +|z|^2\ge |z|^2.$ Then Corollary 8.3
 follows from   $|z| \gtrsim \epsilon  $,  Lemma 6.1.

\bigskip

In the rest of the section we prove Lemma 8.1. We have
$$\aligned &
i\dot \omega = a(\omega , \zeta  )+\langle g ,
    A(\omega ,\zeta
 )\rangle  +
\Cal A(\omega ,\zeta ,g
 )  \\&
i\dot \zeta  -\lambda (\omega )\zeta =d(\omega ,|\zeta |^2 )\zeta
 +   \Cal
B(\omega ,\zeta ,g
 )   \\& i  \dot g -
(\Cal H_\omega +\dot \gamma P_c(\Cal H_\omega )  \sigma _3) g  =
  \widehat{\Cal F}_1
(\omega, \zeta )g  +O_{ loc}(|g|^2) +O(|g|^7). \endaligned \tag 8.3
$$
 We define   recursively, where  $\omega _1=\omega $ and
 (8.4) below for $\ell =1$  coincides
 with the equation for $\omega $ in (8.3),
$$\aligned &     i\dot \omega _{\ell }=
 a_{\ell } (\omega  , \zeta   )+\langle g ,
A_{\ell }(\omega  ,\zeta
 )\rangle  +
\Cal A_{\ell }(\omega  ,\zeta ,g
 )  .  \endaligned \tag 8.4  $$
Let $\delta _2 = \delta _1^2$. We assume:

 {\item {(1)
}}   for $|\omega -\omega _0|\le   e^{2^{-\ell} -2}\delta _2 $ and
$|\zeta |\le \delta _2  $ (the equalities below define
$\widetilde{A}_ \ell$ and $\widetilde{a}_ \ell$)
$$\aligned &  a_\ell (\omega , \zeta )=\sum _{m+n\ge \ell +1}
 a_{\ell mn}(\omega ) \zeta ^m  \overline{\zeta}^n  \, , \,
  \widetilde{a}_ \ell (\omega , \zeta ):=\sum _{m+n\ge \ell +2}
 a_{\ell mn}(\omega ) \zeta ^m  \overline{\zeta}^n \\&  A_\ell (\omega , \zeta )=\sum _{m+n\ge \ell +1}
 A_{\ell mn}(\omega ) \zeta ^m  \overline{\zeta }^n  \, , \,
  \widetilde{A}_ \ell (\omega , \zeta ):=\sum _{m+n\ge \ell +2}
 A_{\ell mn}(\omega ) \zeta ^m  \overline{\zeta}^n;
 \endaligned\tag 8.5 $$

{\item {(2) }}  the following estimates hold for a fixed $C$
$$\aligned & |a _\ell (\omega , \zeta  )|
\le C e^{ (\ell -1)(2-2^{-\ell})} \delta _2 ^{-\ell +1 }    |\zeta
|^{ \ell +1 }\, , \\&   \|A_{\ell} (\omega ,\zeta
 ) \| _{\ell ^{2,2}}\le C e^{ (\ell -1)(2-2^{-\ell})} \delta _2
^{-\ell +1 }   |\zeta  |^{ \ell   }  \, , \quad
   |\Cal A_{\ell} (\omega ,\zeta
 ,g
 ) |\le C \| g \| _{\ell ^{2,-2}}^2.
\endaligned \tag 8.6 $$
These facts hold   for $\ell =1$. We set inductively

$$\aligned & \omega _{\ell +1  }=\omega _{\ell     }+\phi _{\ell +1
}(\omega  , \zeta  ) + \langle g, \Phi _{\ell +1 }(\omega  , \zeta
)\rangle \\& \phi _{\ell +1  }(\omega  , \zeta  )= \sum _{m+n=\ell
+1, m\neq n}
 \frac{a_{\ell mn}(\omega  )
\zeta ^m  \overline{\zeta}^n  }{ (m-n )\lambda (\omega )}
  \\&  \Phi _{\ell +1 }(\omega  , \zeta
)= \sum _{m+n=\ell  }    \zeta ^m \overline{\zeta }^n  R_{\Cal H
_\omega ^*}((n-m )\lambda (\omega )) A _{\ell mn}(\omega  ).
\endaligned \tag 8.7 $$  $\omega _{\ell +1}$ satisfies (8.4) with
$$\aligned    a_{\ell +1} (\omega , \zeta )&=
\widetilde{ a}_\ell (\omega , \zeta ) -\phi _{\ell +1  }'(\omega ,
\zeta  ) (a(\omega , \zeta ) ,  d(\omega , |\zeta |^2)\zeta  ) \, ;
\\    \langle g ,A_{\ell +1} (\omega , \zeta )\rangle  &=
\langle g ,\widetilde{ A}_\ell (\omega , \zeta )\rangle -
\partial _\omega \phi _{\ell +1 }  (\omega , \zeta  )  \langle g ,
    A(\omega ,\zeta
 )\rangle -\\& -\langle
 \left (\dot \gamma P_c(\Cal H_\omega )  \sigma _3+\widehat{\Cal F}_1
(\omega, \zeta )\right )  g  ,
 \Phi _{\ell +1 }(\omega , \zeta
)\rangle; \\    \Cal A_{\ell +1} (\omega , \zeta ,g ) &= \Cal
A_{\ell } (\omega , \zeta ,g )
  - \phi _{\ell +1 }' (\omega , \zeta  )  (
    \Cal A(\omega ,\zeta
 ,g) ,
    \Cal B(\omega ,\zeta
 ,g)  ) -\\& - \langle O_{ loc}(|g|^2) +O(|g|^7),\Phi _{\ell +1 }(\omega , \zeta )\rangle .
\endaligned  \tag 8.8$$
More explicitly,
$$ \aligned    &a_{\ell +1} (\omega , \zeta ) =
\widetilde{ a}_\ell (\omega , \zeta ) -\partial _\omega \phi _{\ell
+1 } (\omega , \zeta  )  a(\omega , \zeta ) -\\& - \partial _\zeta
 \phi _{\ell +1 } (\omega , \zeta  )      d(\omega , |\zeta
|^2)\zeta +
\partial _{\overline{\zeta} } \overline{\phi
_{\ell +1 } (\omega , \zeta  )}     d(\omega , |\zeta
|^2)\overline{\zeta}
\endaligned  \tag 8.9$$
The estimates can be derived in a manner similar to above proofs and
we skip them.   We have $\omega _\ell \to  \varpi $ with the latter
satisfying an equation of the form (8.2). This completes the proof
of Lemma 8.1.

\head \S 9 Wave operators and partial   diagonalization for $  \Cal
H _\omega $
\endhead

\noindent We   set $ \Cal H _0= \sigma_3(H+\omega )$ and write $
\Cal H _\omega = \Cal H _0+  B^\ast (\omega ) A $ with $A(\nu
)=\langle \nu \rangle ^{-\tau}$ with $\tau
>3/2$ and
$B^\ast (\nu ,\omega )$ a    $C^2$ function   in $ \omega $ with
values in the space of $2\times 2$ real valued matrices. For any
$\omega  $ in some compact set $K$  there is a  constant  $c (K)>0$
and $\alpha
>0$ such that
$   \big | e^{\alpha |\nu |}  B^\ast (\nu ,\omega ) \big | \le
c_m(K)(\omega -E_0)
 $ $  \forall   \nu \in \Bbb  Z. $

\proclaim{Lemma 9.1} For $\tau >1$ there exists $C=C(\tau , \omega
)$ such that  for all $z\in \Bbb C \backslash \sigma _e(\Cal
H_\omega )$
$$\| R_{\Cal H_\omega }(z )P_c(\Cal H_\omega   )   \| _{
B (\ell ^{2, \tau }, \ell ^{2,-\tau }) }\le  C   .\tag 9.1$$ There
is a neighborhood $U$ of $\sigma _e(\Cal H_\omega )$ in $\Bbb C$
such that for any $z\in U \backslash \sigma _e(\Cal H_\omega )$
$$\| R_{\Cal H_\omega }(z )   \| _{
B (\ell ^{2, \tau }, \ell ^{2,-\tau }) }\le  C   .\tag 9.2$$

The following limits are well defined
$$\lim _{\epsilon \to 0 ^+}  R_{\Cal H_\omega }(\lambda \pm i\epsilon )
=R^{\pm}_{\Cal H_\omega  }(\lambda  )  \text{ in $C^0(\sigma _e(\Cal
H_\omega ),B(\ell ^{2,\tau } , \ell ^{2,-\tau } ))$} .\tag 9.3$$
\endproclaim
{\it Proof.} First of all (9.2)  implies (9.1). By Lemma 5.7
\cite{CT}  we have
$$ \| \langle x \rangle ^{-\tau} R_{H}( z ,\cdot ,\cdot )
 \langle y \rangle ^{-\tau} \| _{\ell ^2 (\Bbb Z^2)}
  \le  C  \text{ for $z$ close to $[0,4]$}  \tag 1$$
with the following limits well defined

$$\lim _{\epsilon \to 0 ^+}  R_{  H  }(\lambda \pm i\epsilon )
=R^{\pm}_{H  }(\lambda  )  \text{ in $C^0([0,4],B(\ell ^{2,\tau } ,
\ell ^{2,-\tau } ))$} .\tag 2$$ This implies
$$ \| \langle x \rangle ^{-\tau} R_{\Cal H _0}( z ,\cdot ,\cdot )
 \langle y \rangle ^{-\tau} \| _{\ell ^2 (\Bbb Z^2)}
  \le  C  \text{ for $z$ close to $\sigma _e(\Cal H_\omega )$}  \tag 3$$
with the following limits well defined

$$\lim _{\epsilon \to 0 ^+}  R_{ \Cal H _0  }(\lambda \pm i\epsilon )
=R^{\pm}_{\Cal H  _0}(\lambda  )  \text{ in $C^0(\sigma _e(\Cal
H_\omega ),B(\ell ^{2,\tau } , \ell ^{2,-\tau } ))$} .\tag 4$$ For
$A=\langle x \rangle
 ^{-\tau}$ we write
$$AR_{\Cal H_\omega}(z) =
 ( 1+ AR_{\Cal H_0}(z)B^\ast   )^{-1}
  AR_{\Cal H_0}(z) .\tag 5
$$
(3) implies $\| AR_{\Cal H_0}(z)B^\ast \| _{B(\ell ^2,\ell ^2)}\le C
(\omega -E_0)\ll 1  .$ (3)--(4) imply (9.2)--(9.3). This concludes
the proof of Lemma 9.1.
\bigskip

 We have:

 \proclaim{Lemma 9.2} We have:
{\item {(A)}}
 $\Cal H_\omega$ does not have
 resonances at $\pm \omega$ and at   $\pm (4+\omega )$.

 {\item {(B)}} $\sigma _e(\Cal H_\omega )$ does not
 contain eigenvalues.

{\item {(C)}} There are isomorphisms inverses of each other
 $W(\omega )\in B( \ell ^2_c( \Cal  H_{0 })  ,\ell ^2_c( \Cal  H_{\omega })) $ and
$Z(\omega )\in B( \ell ^2_c( \Cal  H_{\omega })  ,  \ell ^2_c( \Cal
H_{0 })  )
  $,
defined as follows:
  \noindent for $u\in \ell ^2_c( \Cal  H_{0 })  $, and $v$ such that $
  \sigma _3v\in \ell ^2_c(  \Cal H_{\omega }
  )  $,
$$\aligned & \langle W u,v\rangle =
\langle  u,v\rangle  +\lim _{\varepsilon \to 0^+ } \frac 1 {2\pi i }
\int _{-\infty}^{+\infty} \langle A R_{\Cal H_0}(\lambda
+i\varepsilon ) u,B R_{\Cal H_\omega ^\ast }(\lambda +i\varepsilon
)v\rangle d\lambda ;
\endaligned$$
for $u\in \ell^2_c( \Cal  H_{\omega }) $,  $v\in  \ell ^2_c( \Cal
H_{0 })  $,
$$\aligned & \langle Z u,v\rangle =
\langle  u,v\rangle  +\lim _{\varepsilon \to 0^+ } \frac 1{2\pi i }
\int _{-\infty}^{+\infty} \langle A R_{\Cal H_\omega }(\lambda
+i\varepsilon ) u,B R_{\Cal H_0}(\lambda +i\varepsilon )v\rangle
d\lambda .
\endaligned$$ Then $ P_c(  H _\omega
) \Cal H_\omega =W \Cal H_0Z.$
  $\| W(\omega )\| _{B(\ell ^2_c( \Cal  H_{0 })  ,\ell ^2_c(  \Cal H_{\omega }   )
)}$ and   $\| Z(\omega )\| _{B(\ell ^2_c(  \Cal H_{\omega }   )  ,
\ell ^2_c( \Cal  H_{0 }) ) }$ are uniformly locally bounded in
$\omega $.
\endproclaim
{\it Proof.} (A) and (B) follow by standard arguments from the fact
that $V_\omega$ is small. (C) follows from  (1)--(4) below.
Specifically we need to show that  there is a fixed  $c>0$ such that
$\forall \, \epsilon \neq 0$
$$\align  & \int \|  \langle x
\rangle ^{-\tau} R_{\Cal  H_{0 }}(\lambda +i\varepsilon)
 u\|^2_{\ell ^2} d\lambda \le
c \| u\| ^2_{\ell ^2}  \text{ for all $u\in \ell ^2_c( \Cal  H_{0 })
$}\tag 1\\& \int \| B R_{\Cal  H_{0 }}(\lambda +i\varepsilon)
u\|^2_{\ell ^2} d\lambda \le c \| u\| ^2_2\text{ for all $u\in \ell
^2_c( \Cal  H_{0 }) $}\tag 2\\& \int \|  B R_{\Cal H_\omega ^\ast
}(\lambda +i\varepsilon ) u\| ^2_{\ell ^2} d\lambda \le c \|
u\|^2_{\ell ^2} \text{ for all $u\in \ell ^2(\Cal H _\omega ^\ast
):=\sigma _3\ell ^2(\Cal H _\omega   )$} \tag 3
\\& \int \|  \langle x \rangle ^{-\tau} R_{\Cal H_\omega }(\lambda
+i\varepsilon)  u\|^2_{\ell ^2} d\lambda \le c \| u\|^2_{\ell ^2}
\text{ for all $ u \in \ell ^2_c(\Cal H _\omega ) $} . \tag 4
\endalign
$$
(1)--(4) are consequences of (9.1) and of inequalities (2)--(3) in
Lemma 9.1.

\bigskip
 We have:

\proclaim{Lemma 9.3} For any  $u\in \ell ^2$ we have
$$\aligned & P_c( \Cal H _{\omega} )u= \lim _{\epsilon \to 0^+}
 \frac 1{2\pi i}  \int _{\sigma _e(\Cal H_\omega )}
  \left [ R_{\Cal H_\omega }(\lambda +i\epsilon
)- R_{\Cal H_\omega }(\lambda -i\epsilon )
 \right ]  ud\lambda   .
\endaligned \tag 1
$$
\endproclaim
{\it Proof.} By $\sigma _e(\Cal H_\omega  )=\sigma _e(\Cal H_0 )= [
\omega , 4+\omega ]\cup [-4-\omega , -\omega ]$ and by the spectral
theorem we have
$$\aligned & P_c( \Cal H _{0} )v= \lim _{\epsilon \to 0^+}
 \frac 1{2\pi i}  \int _{\sigma _e(\Cal H_\omega  )}
  \left [ R_{\Cal H_0 }(\lambda +i\epsilon
)- R_{\Cal H_0 }(\lambda -i\epsilon )
 \right ]  vd\lambda   .
\endaligned \tag 2
$$
To prove (1) in Lemma 9.3 we observe that for $u\in \ell ^2 _d(\Cal
H_\omega )$ both sides of (1) are 0. Hence it is enough to prove (1)
for $u\in \ell ^2 _c(\Cal H_\omega )$. Then $u=W(\omega )v$ with
$v\in \ell ^2 _c(\Cal H_0 )$. Then (1) $u\in \ell ^2 _c(\Cal
H_\omega )$ can be rewritten as
$$\aligned & u= W(\omega )v= W(\omega )\lim _{\epsilon \to 0^+}
 \frac 1{2\pi i}  \int _{\sigma _e(\Cal H_\omega )}
  \left [ R_{\Cal H_0 }(\lambda +i\epsilon
)- R_{\Cal H_0 }(\lambda -i\epsilon )
 \right ]  vd\lambda    \\& =
\lim _{\epsilon \to 0^+}
 \frac 1{2\pi i}  \int _{\sigma _e(\Cal H_\omega )}
  \left [ R_{\Cal H_\omega }(\lambda +i\epsilon
)- R_{\Cal H_\omega }(\lambda -i\epsilon )
 \right ]  W(\omega )vd\lambda  .
\endaligned
$$
  This concludes Lemma
9.3.

\bigskip
Finally, we obtain the limiting absorption principle:

\proclaim{Lemma 9.4} For any  $u\in \ell ^{2,\tau } $ with $\tau
>3/2$
we have
$$\aligned & P_c( \Cal H _{\omega} )u=
 \frac 1{2\pi i}  \int _{\sigma _e(\Cal  H_\omega )}
  \left [ R_{\Cal  H_\omega }^+(\lambda
)- R_{\Cal  H_\omega }^-(\lambda   )
 \right ]  ud\lambda   .
\endaligned
$$
\endproclaim
{\it Proof.} For $u\in \ell ^{2,\tau } $   the $\epsilon \to 0^+$
limit in (1)  Lemma 9.4 converges  in $\ell ^{2,-\tau } $ to the
integral in Lemma 9.5.

\bigskip
 \head \S 10 Dispersive theory for $\Cal H _\omega $:
 proof of Lemma 3.1 \endhead
Lemma 9.2 implies $ \| P_c(\Cal H_\omega )e^{it\Cal H_\omega } \|
_{B(\ell ^2  , \ell ^{2})}   \le C $ for a fixed $C>0$ and yields
Lemma 3.1 for $p=2$. By interpolation the rest of Lemma 3.1 will be
a consequence of  case $p=1$. Lemma 9.4 implies for $u\in \Cal S
(\Bbb Z)$
$$\aligned & P_c( \Cal H _{\omega} )e^{it\Cal H_\omega }u=
 \frac 1{2\pi i}  \int _{\sigma _e(\Cal  H_\omega )}
  e^{i\lambda t}\left [ R_{\Cal  H_\omega }^+(\lambda
)- R_{\Cal  H_\omega }^-(\lambda   )
 \right ]  ud\lambda  .
\endaligned
$$
We expand $ R_{\Cal  H_\omega }^\pm (\lambda )=\sum _{j=0}^\infty
(-1)^j  \left ( R_{\Cal  H_0 }^\pm (\lambda ) V_\omega \right
)^jR_{\Cal  H_0 }^\pm (\lambda ).$  Correspondingly
$$\aligned &P_c( \Cal H _{\omega} )e^{it\Cal H_\omega }u
=\sum _{j=1}^\infty (-1)^j \frac 1{2\pi i}  \int _{\sigma _e(\Cal
H_\omega )}
  e^{i\lambda t}  \left ( R_{\Cal  H_0 }^+(\lambda )
V_\omega \right )^jR_{\Cal  H_0 }^+ (\lambda )u d\lambda \\& - \sum
_{j=1}^\infty (-1)^j   \frac 1{2\pi i}  \int _{\sigma _e(\Cal
H_\omega )}
  e^{i\lambda t}  \left ( R_{\Cal  H_0 }^- (\lambda )
V_\omega \right )^jR_{\Cal  H_0 }^- (\lambda )ud\lambda +P_c( \Cal H
_{0} )e^{it\Cal H_0 }u .
\endaligned
$$
We have $  \| P_c( \Cal H _{0} )e^{it\Cal H_0 } \| _{B(\ell ^1  ,
\ell ^{\infty})}   \le C\langle t \rangle ^{-\frac 13} $ by
\cite{CT,PS}. We consider now a generic term in the above
summations. By  $\sigma _e(\Cal H_\omega  )=  [ \omega , 4+\omega
]\cup [-4-\omega , -\omega ]$ it is not restrictive to focus on what
follows, for $K(E)=\left  ( R_{\Cal  H_0 }^+(E+\omega  ) V_\omega
\right )^j $ with  $j\ge 0$,
$$\aligned
& Tu=\int _{0}^4
  e^{iE t}    R_{\Cal  H_0 }^+(E+\omega  )
V_\omega K(E)R_{\Cal H_0 }^+ (E+\omega )u dE.
\endaligned $$
Let $P^0_+=\text{diag}(1,0)$ and $P^0_-=\text{diag}(0,1)$. Then $
T=P^0_+TP^0_++P^0_-TP^0_++P^0_+TP^0_-+P^0_-TP^0_ -.$ We consider
separately these four operators. We focus on $T_{++}=P^0_+TP^0_+$.
We have
$$\aligned & T_{++}(\mu , \nu )=\sum _{(\nu ', \mu ')\in \Bbb Z^2}
\int _{0}^4
  e^{i\lambda  t}   P^0_+ R_{   H  }^+(\lambda  , \nu , \nu ' )
V_\omega (\nu ')K(\lambda ,  \nu ' , \mu ' )R_{  H  }^+ (\lambda ,
\mu ' , \mu  ) P^0_+ d\lambda .
\endaligned $$
Recall now, see for instance Lemma 5.9 \cite{CT}, that
 $$ \aligned & R_{   H  }^+(\lambda , \nu , \nu ' )= -
 \frac{ {f}_+(\nu ,
\theta )\, {f}_- (\nu ',  \theta ) }{W(  \theta )} \quad \text{for
$\nu \ge \nu  '$ ,}
\\ & R_{   H  }^+(\lambda , \nu , \nu ' )=-\frac{ {f}_-(\nu , \theta )
 {f}_+ (\nu ', \theta )  }{W(  \theta )} \quad \text{for $\nu <\nu '
  $ }
\endaligned \tag 10.1
$$
where $\theta \in [0,\pi]$ is such that $\lambda =2(1-\cos \theta
)$,   where ${f}_\pm  (\nu  ,  \theta )$ are the Jost functions
satisfying
$$H{f}_{  \pm } (\nu , \theta )=\lambda {f}_{  \pm }
 (\nu , \theta ) \text{ with }
\lim _{\nu \to \pm \infty}\left [ {f}_{  \pm } (\nu , \theta )
 - e^{\mp
i\nu  \theta }\right  ]=0  $$ and where $W(  \theta )= {f}_+  (\nu
+1 , \theta ){f}_- (\nu  ,  \theta )-{f}_+ (\nu  ,  \theta ) {f}_-
(\nu +1 , \theta )$ is the   Wronskian . We recall that $W(\theta )$
is $C^\infty $ in $\Bbb R/2\pi \Bbb Z$ and that $W(\theta )\neq 0$
for all $\theta$, see Lemmas 5.3 and 5.5 in \cite{CT}. We split the
above sum in various terms
$$T_{++}(\mu , \nu )=\sum _{\nu  '\le \nu , \mu   \le  \mu ' }\cdots +
\sum _{\nu  '> \nu , \mu   \le  \mu ' }\cdots + \sum _{\nu  '\le \nu
, \mu   >  \mu ' }\cdots + \sum _{\nu  '> \nu , \mu   > \mu '
}\cdots . \tag 10.2
$$
Let ${m}_{  \pm } (\nu ,\theta )= e^{\pm i \nu \theta
  } {f}_{  \pm } (\nu , \theta )
  .$ Then the first term in the rhs in (10.1) is
 $$\aligned & \sum _{\nu  '\le \nu , \mu   \le  \mu ' }
 \int _{ 0 }^\pi e^{it  (2-2\cos \theta )+i(\nu -\mu )\theta } {m}_- (\nu  ,  \theta )
 A(\theta ,  \nu ' , \mu ' )
{m}_+ (\mu,  \theta )  \sin \theta d\theta \endaligned \tag  10.3 $$
with $A(\theta ,  \nu ' , \mu ' )= {f}_+ (\nu ' ,  \theta )V_\omega
(\nu ')K(\lambda ,  \nu ' , \mu ' ) {f}_- (\mu ' ,  \theta )$. We
have
$$\aligned & \sup _{\nu  '\le \nu , \mu   \le  \mu ' }
\left ( \langle \nu '\rangle ^2 \langle \mu '\rangle ^2\| {m}_- (\nu
, \theta )
 A(\theta ,  \nu ' , \mu ' )
{m}_+ (\mu,  \theta ) \| _{W^{1,1}(0,\pi )}\right ) \le C^j (\omega
-E_0)^{j+1} .
\endaligned $$
Then by stationary phase $|(10.3)|\le C^j (\omega -E_0)^{j+1}
\langle t \rangle ^{-\frac{1}{3}}$.  The other terms in the rhs of
(10.2)  can be bounded using
$$ {f}_{  \mp } (\nu ,\theta )=
\frac{1}{ T  (\theta )}\overline{ {f}_{  \pm } (\nu ,\theta )}
+\frac{R_{ \pm} (\theta)}{T  (\theta  )}  {f}_{  \pm } (\nu ,\theta
)   $$ where transmission and reflection coefficients are defined
using the Wronskians

$$T  (\theta )= \frac {\pm 2 i\sin \theta } {[ f_{ \mp
} (\theta ),
 f_{ \pm } (\theta )]  } \, , \quad R_{  \pm } (\theta ) = -\frac {[
 {f_{ \mp } ( \theta  )},  \overline{f_{  \pm }} (\theta)] }
{[ f_{ \mp } ( \theta ),  f_{  \pm } (\theta )] },$$ see \cite{CT}.
Finally, the remaining terms $P^0_-TP^0_+$, $P^0_+TP^0_-$ and
$P^0_-TP^0_ - $ can be bounded similarly. So we have $\| T\|
_{B(\ell ^1  , \ell ^{\infty})}\le C^j (\omega -E_0)^{j+1} \langle t
\rangle ^{-\frac{1}{3}}$. By summing up on $j$ and for $C  (\omega
-E_0)<1$ we obtain $  \| P_c( \Cal H _{\omega} )e^{it\Cal H_\omega }
\| _{B(\ell ^1  , \ell ^{\infty})}   \le C\langle t \rangle ^{-\frac
13} $.

\head \S 11 The projections $P_\pm (\omega )$\endhead

Let $u\in \Cal S (\Bbb Z)$. We set

$$\aligned & P_\pm (  {\omega} )u=
 \frac 1{2\pi i}  \int _{\sigma _e(\Cal  H_\omega )\cap \Bbb R_\pm }
  \left [ R_{\Cal  H_\omega }^+(\lambda
)- R_{\Cal  H_\omega }^-(\lambda   )
 \right ]  ud\lambda   \\& = \lim _{\epsilon \to 0^+}
 \frac 1{2\pi i}  \int _{\sigma _e(\Cal H_\omega )\cap \Bbb R_\pm }
  \left [ R_{\Cal H_\omega }(\lambda +i\epsilon
)- R_{\Cal H_\omega }(\lambda -i\epsilon )
 \right ]  ud\lambda .
\endaligned
$$
Notice that we have $P_\pm (  {\omega} ) =Z(\omega ) P_\pm ^0   P_c
(H) W({\omega} )$.  Hence $P_\pm (  {\omega} )$ extend into
projections in $\ell ^2$. During the course of the proof of  Lemma
5.1 we used the following fact, which we prove:

\proclaim{Lemma 11.1} For any pair $s_1,s_2\in \Bbb R$ there is $
c_{s_1,s_2} (\omega )$ upper semicontinuous in $\omega$  such that
for  $j=0,1$

$$   \|   P_c (\omega
) \sigma _3- (P_+ (\omega )-P_- (\omega ) \| _{B (\ell ^{2,s_1}
,\ell  ^{2,s_2}  ) } \le c_{s_1,s_2} (\omega )<\infty . \tag 1$$

\endproclaim
{\it Proof.} For this proof we set $\Cal H=\Cal H_\omega$, $\Cal
H_0=\sigma _3 (H + \omega )$,
  $R_0(z)= (\Cal H_0-z)^{-1}$ and $R(z)= (\Cal H-z)^{-1}$.
To prove (1) it is enough to write $P_c=P_++P_-$ and to prove
$\|\left [ P_\pm \sigma _3 \mp P_\pm \right ] g\|  _{\ell ^ {2, M} }
\le c\| g\| _{\ell^ {2, -N} } .$ It is not restrictive to consider
only $P_+$. Setting $H=H_0+V$, we write
$$\aligned &\sum_\pm \pm R(\lambda \pm  i\epsilon )
= \sum _\pm \pm (1+ R_0(\lambda \pm  i\epsilon ) V )^{-1}
R_0(\lambda \pm  i\epsilon ) .
\endaligned \tag 2
$$
By elementary computation
$$  R_0(\lambda \pm  i\epsilon )
\sigma _3 =R_0(\lambda \pm  i\epsilon ) -2 (H +\omega +\lambda \pm i
\epsilon )^{-1} \text{diag}(0,1).$$ Therefore  $$ \text{rhs}  (2)
\sigma _3 = \text{rhs} \,(4)+2\sum _\pm \pm (1+ R_0(\lambda \pm
i\epsilon ) V )^{-1} \text{diag}(0,1)
  R_{H}   (- \omega -\lambda \mp i
\epsilon ) .
$$
Hence we are reduced to show that $ Ku=$
$$
\lim _{\epsilon \to 0^+} \sum _\pm \pm \int _{\Bbb R_+\cap \sigma _c
(\Cal H_0)} (1+ R_0(\lambda \pm  i\epsilon ) V )^{-1}
\text{diag}(0,1)
 R_{H}   (- \omega -\lambda \mp i
\epsilon )
 u d\lambda
$$
defines an operator such that  for some fixed $c$  $$\| Ku \| _{\ell
^ {2, M} } \le c \| u \| _{\ell ^ {2, -N} } \tag 3$$
 We expand $(1+R_0V )^{-1}=\sum
_{j=0}^{\infty } \left [ - R_0V \right ] ^j  $ and we consider the
corresponding decomposition $ K=\sum _{j=0}^{\infty }K_j^0 . $ We
have $K_0^0= 0$ since for any $u\in \ell ^ {2 }$ we have
$$  \lim _{\epsilon \to 0^+}  \int
_{\Bbb R_+\cap \sigma _c(\Cal H_0)}
 \sum _\pm \pm (H +\omega +\lambda \pm i
\epsilon )^{-1} \text{diag}(0,1)ud\lambda =0.$$ We next consider
$K_j^0$  for $j>0$ and prove $$\| K_j^0u \| _{\ell ^ {2,M }} \le
c^j(\omega -E_0)^j \| u \| _{\ell ^ {2,  -N} } \text{ for a fixed
$c$.} \tag 4$$ We have that for $\gamma $ a closed path around
$[\omega , 4+\omega ]$ we have
$$K_j^0u =
\int _{\gamma }  R_0(z ) V(R_0(z ) V )^{j-1}  \text{diag}(0,1)
 R_{H}   (- \omega -z )
 u dz.
$$
We have $\| V(R_0(z ) V )^{j-1}\| _{B(\ell ^ {2,-N }, \ell ^ {2,M
})}\le \widetilde{c}^j(\omega -E_0)^j$. For  $z\in \gamma$ we have
$$ | R_0(z,\nu , \mu )|+ |R_{H}   (- \omega -z ,\nu , \mu )|
\le  \beta e^{-\alpha |\nu -\mu |}$$ for some fixed $\alpha >0$ and
$\beta
>0$. Then
$$\aligned & \| K_j^0\| _{B(\ell ^ {2,-N }, \ell ^ {2,M
})} \\& \le   \beta ^2 \| e^{-\alpha |\cdot  |}*\| _{B(\ell ^ {2,-N
}, \ell ^ {2,-N })} \| V(R_0(z ) V )^{j-1}\| _{B(\ell ^ {2,-N },
\ell ^ {2,M })} \| e^{-\alpha |\cdot  |}*\| _{B(\ell ^ {2,-N }, \ell
^ {2,M })} \\&  \le   {c}^j(\omega -E_0)^j.\endaligned $$ This
yields Lemma 11.1.

\head \S Appendix A: existence of $H$ satisfying  (H1)--(H3)\endhead

Lemmas A.1 and A.2 below, together    prove the existence of $H$
satisfying (H1)--(H3). Lemma A.1  uses standard perturbation
arguments.  Lemma A.2 proves that, in some sense, 0 average
potentials generically satisfy the second condition in (1) below. We
recall that the resolvent $R_{-\Delta}(z)$ for $ z\in \Bbb C
\backslash [0,4]$ has kernel
$$
R_{-\Delta}(\mu , \nu , z)=\frac{1}{2i\sin \theta }e^{-i \theta
|\mu -\nu |}, \ \ \mu ,\nu  \in \Bbb Z,
$$
with $\theta$ a solution     to $    2(1-\cos \theta ) =z
 $ with $\Im \theta \le 0$. We will consider $q(\mu )\in S(\Bbb
 Z)$.  \proclaim{Lemma A.1} Suppose
 that for some $a_0>0$ and $C_0$ we have $|q(\mu )|\le
C_0e^{-a_0|\mu |} $ and that the following two conditions are
satisfied:
$$\sum _{\mu \in \Bbb Z}q(\mu )=0 \, ,  \,
\sum _{\mu ,\nu \in \Bbb Z}|\mu -\nu |q(\mu )q(\nu )<0 .\tag 1$$
Then there is $\varepsilon _0>0$ such that for any $0<\varepsilon
<\varepsilon _0$ the operator $H=-\Delta +\varepsilon q$ satisfies
hypotheses (H1)--(H3).
\endproclaim
{\it Proof.} For $z\not \in \sigma (H)$ we have
 $  R_{H}(z)= (1+R_{-\Delta}(z)\varepsilon q)^{-1}R_{-\Delta}(z).$
So there are two eigenvalues, one in $(-\infty ,0)$ and the other in
$(4,\infty )$, exactly if $(1+R_{-\Delta}(z)\varepsilon q)^{-1}$ is
singular in two such points. By Fredholm theory the singular points
occur in correspondence to values of $z$ such that $\ker
(1+R_{-\Delta}(z)\varepsilon q)\neq 0$. If we set $b(\nu )=\sqrt{
|q(\nu )|}$ and $a(\nu )=b(\nu )\text{sign $q(\nu )$}$ we have
$q(\nu )=a(\nu )b(\nu ).$ By standard arguments the map $u\to v=bu$
establishes an isomorphism $\ker (1+R_{-\Delta}(z)\varepsilon q)\to
\ker (1+bR_{-\Delta}(z)\varepsilon a) $ with inverse $u=-R_{-\Delta
}(z)av$. So $\ker (1+R_{-\Delta}(z)\varepsilon q)\neq 0$ exactly if
$\ker (1+bR_{-\Delta}(z)\varepsilon a)\neq 0.$ We split the analysis
in two parts. We first look at $z$ near 0. We will later look at $z$
near $4$. For $z$ near 0 write
$$R_{-\Delta}(\mu , \nu , z)=\frac{1}{2i\sin \theta } -D_0(\mu , \nu , z)
\, , \quad D_0(\mu , \nu , z)=\frac{1  -e^{-i \theta |\mu -\nu
|}}{2i\sin \theta } .$$ Now we have
$$\aligned & 1+bR_{-\Delta}(z)\varepsilon a=
1+\varepsilon \frac{b  \langle \cdot , a\rangle}{2i\sin \theta} -
\varepsilon bD_0(z)a =\\& = \frac{1}{2i\sin \theta }(1-\varepsilon
bD_0(z)a) \left ( 2i\sin \theta  + \varepsilon (1-\varepsilon
bD_0(z)a)^{-1}b \langle \cdot , a\rangle \right ) .
\endaligned $$
We are reduced at looking at the kernel of the third factor in the
last line. We will show that the following equation admits a
solution $\theta =it$ with $t<0$ close to 0. The singularity we are
searching corresponds to solutions of
$$ \aligned & 2i\sin \theta  +\langle (1-\varepsilon
bD_0(z)a)^{-1}b , a\rangle =    0= 2i\sin \theta  +\varepsilon  \sum
_{n=1}^{ \infty}  {\varepsilon ^n} \langle (bD_0(z) a)^nb, a\rangle
\\& =2i\sin \theta  + \frac{\varepsilon ^2}{2} \sum _{\mu ,\nu \in \Bbb
Z}|\mu -\nu |q(\mu )q(\nu ) +O(\varepsilon ^3)+O(\theta ^2).
\endaligned
$$
  By
the implicit function theorem we obtain a unique solution
$$\theta
=\theta ( \varepsilon )  = i\frac{\varepsilon ^2}{4} \sum _{\mu ,\nu
\in \Bbb Z}|\mu -\nu |q(\mu )q(\nu )  +O(\varepsilon ^2).$$ For
$\varepsilon ^2
>0$ we have $\Im \theta ( \varepsilon )<0$. This yields for $\varepsilon
^2>0$ an eigenvalue $z(\varepsilon )=2(1-\cos (\theta ( \varepsilon
)))$ of $H$ with $z(\varepsilon )$ close to 0. Necessarily $
z(\varepsilon )<0$ by Lemma 5.3 \cite{CT} and by selfadjointness of
$H$. So we obtain $\theta =it$ with $t<0$ and $t=O(\varepsilon ^2)$.
This argument proves the existence of an eigenvalue also near 4,
since if $(-\Delta +q) u=\lambda u$ then $v(\nu )=(-1)^\nu u(\nu )$
satisfies $(-\Delta -q) v=(4-\lambda )v$.   For $\varepsilon _0$
small, there are no other eigenvalues. We need to show that $0$ and
$4$ are not resonances. We focus on 0.  0 for some $\varepsilon $ is
a resonance exactly if the Wronskian $W( \varepsilon,  \theta ) $ is
$W( \varepsilon,  0 ) =0$ (see \S 10, we have added the parameter
$\varepsilon$). By (10.1) this happens exactly if
$$\aligned & bR_H(z)a=
(1+\varepsilon bR_{-\Delta}(z) a)^{-1}bR_{-\Delta}(z)a=\\&
  \left ( 2i\sin
\theta  + \varepsilon (1-\varepsilon bD_0(z)a)^{-1}b \langle \cdot ,
a\rangle \right ) ^{-1} (1-\varepsilon bD_0(z)a) ^{-1}  \left (
{2i\sin \theta } bR_{-\Delta}(z)a \right )
 \endaligned \tag  1
$$ is singular at $z=0$. But the first factor in rhs(1) is not
singular at $z=\theta =0$. The second factor is not singular. The
third factor has kernel

$$ k(\mu , \nu )= b(\mu ) a(\nu ) -   b(\mu )(1  -e^{-i \theta |\mu -\nu
|}) a(\nu )  $$ which is obviously not singular at $\theta =0$. This
means that 0 is not a resonance. The same argument proves also the
statement for 4, thanks to the   transformation $v(\nu )=(-1)^\nu
u(\nu )$. This yields Lemma A.1.
\bigskip

Here we recall that  \cite{KKK} observe that if   $q(\mu )=A \delta
(\mu - \mu _1)+B\delta (\mu - \mu _2)$ with $(A,B)\neq (0,0)$ then
$H=-\Delta +   q$ has at least one eigenvalue. In the case of small
potentials we can generalize this observation. If $ \sum _{\mu \in
\Bbb Z}q(\mu )\neq 0$ then proceeding as in Lemma A.1 it is easy to
show that $H=-\Delta +\varepsilon q$ has an eigenvalue for
$\varepsilon \neq 0$ small.  For $ \sum _{\mu \in \Bbb Z}q(\mu )= 0$
we can generally apply Lemma A.1 thanks to the following lemma:

 \proclaim{Lemma A.2}
Suppose that   $q(\mu )\in \Cal S(\Bbb
 Z)$  satisfies $ \sum _{\mu \in \Bbb Z}q(\mu )=0$. Then
$$   -\frac{1}{2}\sum _{\mu ,\nu \in \Bbb Z}|\mu -\nu |q(\mu )q(\nu
) =\lim _{z \to  0^-} \langle R _{-\Delta} (z)q,q\rangle .$$
\endproclaim
Notice that for $z<0$ we have $\langle R _{-\Delta} (z)q,q\rangle
>0$, so  the above limit is generically positive.
To prove Lemma A.2 use
$$
R_{-\Delta}(\mu , \nu , z)=\frac{1}{2i\sin \theta }e^{-i \theta |\mu
-\nu |}= \frac{1}{2i\sin \theta } - \frac{\theta |\nu -\mu |}{2 \sin
\theta}+O(\theta  ),
$$
We have for $|\cdot |\ast q (\mu )= \sum _\nu |\mu - \nu | q(\nu )$
$$\aligned & \langle R _{-\Delta} (z)q,q\rangle  =
\frac{\sum _\nu q(\nu )  }{2i\sin \theta }-\frac{1}{2}\langle |\cdot
|\ast q,q \rangle +O(\theta )\to -\frac{1}{2}\langle |\cdot |\ast
q,q \rangle .
\endaligned$$

\enddocument